\documentclass[12pt,a4paper,oneside]{amsart}
\usepackage[cp1251]{inputenc}
\usepackage{amsfonts, amsmath, amssymb, amsthm, amscd, hyperref}
\usepackage[T2A]{fontenc}
\usepackage[cp1251]{inputenc}
\usepackage[english]{babel}
\usepackage[dvips]{graphicx}
\usepackage{subfig}
\usepackage{psfrag}
\usepackage{epstopdf}
\usepackage{cmap}
\usepackage[english]{babel}
\usepackage{amssymb}
\usepackage{amsmath}
\newtheorem{lemma}{Lemma}
\newtheorem{definition}{Definition}
\newtheorem{remark}{Remark}
\newtheorem{fact}{Theorem}
\newtheorem{corollary}{Corollary}
\newtheorem{conjecture}{Conjecture}

\newcommand{\la}{\lambda}

\newcommand{\e}{\varepsilon}

\def\frl{\forall}%
\newcommand{\men}{\leqslant}
\newcommand{\bol}{\geqslant}
\newcommand{\bra}{\langle}
\newcommand{\ket}{\rangle}
\newcommand{\B}{B}
\newcommand{\R}{\mathbb{R}}

\newcommand{\mgl}[2]{\rho_{#1}\left( #2 \right)}
\newcommand{\norm}[1]{\left\| #1 \right\|}

\def\vn{\mathop{\rm int}}
\def\co{\mathop{\rm co}}
\def\Lin{\mathop{\rm Lin}}

\def\dim{\mathop{\rm dim}}

\newcommand{\mglx}[1]{\rho_{X}\!\!\:\left( #1 \right)}\!%
\newcommand{\psif}[1]{\psi\!\left( #1 \right)}%
\newcommand{\mcox}[1]{\delta_{X}\!\left( #1 \right)}%
\newcommand{\lamx}[1]{\lambda^{-}_{X}\!\left( #1 \right)}%
\newcommand{\lapx}[1]{\lambda^{+}_{X}\!\left( #1 \right)}%
\newcommand{\mgli}[2]{\rho_{#1}^{-1}\!\left( #2 \right)}%
\newcommand{\omxi}[1]{\omega_{X}^{-1}\!\left( #1 \right)}%

\newcommand{\cmgl}{\sqrt{2} - 1}%
\newcommand{\omal}[1]{\textit{o}\!\left( #1 \right)}%
\newcommand{\Ker}{\operatorname{Ker}}%
\newcommand{\lin}[1]{\operatorname{Lin}\left\{#1\right\}}%

\newcommand{\prp}{\urcorner}%
\newcommand{\prooff}{\noindent {\bf Proof.}\\}%
\newcommand{\bbox}{\par\noindent\ensuremath{\Box}\par\noindent}%
\newenvironment{prf}
{\prooff}
{\bbox}

\def\B{\mathfrak{B}}%
\def\SS{\partial\B_1(o)}%
\def\SSS{\partial\B_1^*(o)}%
\def\BB{\B_1(o)}%
\newcommand{\reff}[1]{(\ref{#1})}%

\graphicspath{{figure/}}
\newtheorem{nfact}{Theorem}

\def\dd{\ensuremath{\delta}}%
\begin{document}

\title{Hypomonotonicity of the normal cone and proximal smoothness}

\author{G.M. Ivanov}

\address{Department of Higher Mathematics, Moscow Institute of Physics and Technology,  Institutskii pereulok 9, Dolgoprudny, Moscow
region, 141700, Russia}
\address{
National Research University Higher School of Economics,
School of Applied Mathematics and Information Science,
Bolshoi Trekhsvyatitelskiy~3, Moscow, 109028, Russia}
 \email{grimivanov@gmail.com}
\begin{abstract}
In this paper we study the properties of the normal cone to the proximally smooth set.
We give the complete characterization of a proximally smooth set 
through the monotony properties of its normal cone in an arbitrary uniformly convex and uniformly smooth 
Banach space. We give the exact bounds for right-hand side in the monotonicity inequality for normal cone 
in terms of the moduli of smoothness and convexity of a Banach space. 
\end{abstract}

\maketitle

\section{Introduction}
Let $X$ be a real Banach space. For a set $A \subset X$ by $  \partial A, \,\vn A$ 
we denote the boundary and the  interior of  $A.$
We use $\bra p,x \ket$ to denote the value of  functional $p \in X^*$ at   vector $x \in X.$
For $R>0$ and $c \in X$ we denote by $\B_R(c)$ the closed ball with center $c$ and radius $R,$
by $\B_R^{*}(c)$ we denote the ball in the conjugate space.
The distance from a point $x \in X$ to a set $A \subset X$ is defined as 
$$
\rho (x, A) = \inf\limits_{a \in A} \norm{x - a}.
$$
The {\it metric projection} of a point $x $ onto a set $A$ is defined as any element of the set 
$$
	P_A(x) = \{a \in A: \norm{a-x} = \rho(x, A)\}.
$$

Clarke, Stern and Wolenski \cite{Clark} introduced and studied the {\it proximally smooth} sets in
a Hilbert space $H.$ 
\begin{definition}
A set $A \subset X$ is said to be proximally smooth 	with constant $R$ if the distance function
 $x \rightarrow \rho(x, A)$ is continuously differentiable on  set 
$$U(R,A)=\left\{ x \in X: 0 < \rho(x, A) < R \right\}.$$ 
\end{definition}
We denote by $\Omega_{PS}(R)$ the set of all closed proximally smooth  sets
with constant $R$ in $X.$ 

Properties of proximally smooth sets an a Banach space and the relation between such sets and akin 
classes of sets were investigated in \cite{Clark}-\cite{IvBal}. 
Such sets are usually called  weakly convex sets. Some notes about the history of this problem 
one could find in paper \cite{Pol_Rock_Thi1}. Perhaps, the most complete description of the proximally smooth sets 
and their properties in a Hilbert space can be found in  \cite{IvMonEng}.

\begin{definition} \label{def_nconus}
The {\emph{normal cone}} to a set $A \subset X$ at
a point $a_0\in A$ is defined as follows
$$
N(a_0, A) = \Bigl\{ p \in X^*: \frl \e > 0 \; \exists\; \delta
>0: 
\quad  \frl a \in A\cap \B_{\delta}(a_0)\quad
\bra p, a - a_0\ket \le \e \norm{a - a_0} \Bigr\}.
$$
\end{definition}
\begin{remark}\label{remark_def_N}
It is easily seen that normal cone $N(a_0, A)$ does not change when the norm of the space $X$ changes in an equivalent way.
\end{remark} 

The problem of characterization of proximally smooth sets in Hilbert spaces through the monotony properties of the normal cone turns out
to be very important for applications. Thus, Theorem 1.9.1 of monograph \cite{IvMonEng} may be reformulated in the following way (See also \cite{Pol_Rock_Thi1}, Cor. 2.2).
\begin{nfact} \label{th_main_hilbert} 
	Let $A$ be a closed set in a Hilbert space $H$ and $R>0.$ The following conditions are equivalent
	\begin{enumerate}
		\item the set $A $ is proximally smooth with constant $R > 0;$
		\item  for any vectors $x_1,x_2 \in A,$ $p_1 \in N(x_1, A),$ $p_2 \in N(x_2, A)$ such that
		$\norm{p_1}=\norm{p_2} = 1,$ the following inequality holds
		\begin{equation}
			\bra p_2 - p_1, x_2 - x_1 \ket \bol -\frac{\norm{x_2 - x_1}^2}{R}.
			\label{eq:proxreq_cond}
		\end{equation}
	\end{enumerate}
\end{nfact}

The equivalent in a Hilbert space properties of the normal cone were studied earlier in  \cite{Clark}, \cite{proxreg_rockafellar}. 
Note that the fact that  inequality  \reff{eq:proxreq_cond} holds for proximal smooth sets in Hilbert spaces  was the starting point for the construction of the proximal subgradient,  which is an important construction in optimization nowadays. 

Afterwards, of course, attempts were made to extend this result to Banach spaces. Here we must note the works of Thibault and his colleagues, for example,  \cite{Thibault2}, \cite{t1}. As a rule, in these works other approaches to defining the normal cone are examined, and then the coincidence of different normal cones for proximally smooth sets in special Banach spaces  (usually, uniformly smooth and uniformly convex) is proved. 

Perhaps, the strongest results in this domain may be formulated in the following two theorems. 

\begin{nfact}
	Let  $X$ be a uniformly convex and uniformly smooth Banach space.
	Let $\mglx{\tau} \asymp \tau^2$ as $\tau \to 0.$
	Then a proximally smooth set $A\subset X$ with constant $r >0$ satisfies condition 2) of Theorem \ref{th_main_hilbert} for some constant $R>0.$
\end{nfact}
\begin{nfact}
	Let the convexity and smoothness moduli be of power order at zero in a Banach space $X$. 
	Let $\mcox{\e} \asymp \e^2$ as $\e \to 0.$
Then, if a set $A$ satisfies condition 2) of Theorem \ref{th_main_hilbert}, then it is proximally smooth with some constant $ r > 0.$ 
\end{nfact}

In this paper, using the geometrical properties of the sets from the class $\Omega_{PS}(R),$  the connection between the smoothness of a  set and the properties of its normal cone is studied.

\section{Main results}
In the sequel we shall need some additional notation.

By definition, put $J_1(x) = \{p \in \SSS :\, \bra p, x \ket\ = \norm{x}\}.$

Let
$$
\delta_X(\e) = \inf \left\{ 1 - \frac{\|x + y\|}{2}:\ x,y\in\B_{1}(0),\ \|x -y\| \ge \e\right\}
$$
and 
$$
\mglx{\tau} = \sup \left\{\frac{\|x + y\|}{2} + \frac{\|x - y\|}{2} - 1 : \, \|x\| = 1, \|y\| = \tau \right\}.
$$
The functions $\delta_X(\cdot): [0,2] \to [0, 1]$ and $\mglx{\cdot}: \R^+ \to \R^+$
 are  referred to as the moduli of convexity and smoothness of $X$ respectively.

Normed space $X$ is called {\emph{uniformly
convex}}, if  $\mcox{\e} > 0$ for all $\e \in (0, 2].$ 
Normed space $X$ is called {\emph{uniformly
smooth}}, if 
$\lim\limits_{\tau \to +0}{\frac{\mglx{\tau}}{\tau}} = 0.$

Let $f$ and $g$  be two non-negative functions, each one defined on a segment $[0, \e].$
We shall consider  $f$ and $g$ as {\it equivalent at zero,} denoted by $f(t) \asymp g(t)$ as $t \to 0,$ 
if there exist positive constants 
$a,b,c,d,e$ such that $a f(bt) \men g(t) \men c f(dt)$ for $t \in [0,e].$

First we generalize condition 2) of Theorem \ref{th_main_hilbert}.

\begin{definition}
Let a function $\psi:[0,+\infty)\to [0,+\infty)$ be given.
	A set $A\subset X$ satisfies the $\psi$-hypomonotonity condition of the normal cone with constant $R>0$ if for some $\e > 0 $ and for any $x_1, x_2 \in A,$ $p_1 \in N(x_1,\, A),$ $p_2 \in N(x_2,\, A),$ 
	$\norm{p_1} = \norm{p_2} = 1$ such that $\norm{x_1 - x_2} \men \e,$
	the inequality
	\begin{equation}
	\label{cone_to_prox_hypo1}
	\bra p_2 - p_1, x_2 - x_1 \ket \bol - R \psif{\frac{\norm{x_2 - x_1}}{R}}
	\end{equation}
	holds.
	Through $\Omega_{N}^{\psi}(R)$ we denote the class of all closed sets 
	$A\subset X$ that satisfy the $\psi$-hypomonotonity condition with constant $R>0$.
\end{definition}

Naturally, studying the classes of sets that satisfy the
$\psi$-hypomonotonity condition, we suppose some restrictions on the function $\psi.$
\begin{definition}
Through $\mathfrak{M}$ denote the class of convex and Lipschitz functions 
$\psi: \; \R^+ \to \R^+$ such that $\psi(0) = 0.$ 
\end{definition}
We say that  the function  $N : [0, +\infty) \to [0, +\infty)$ such that
$N(0) = 0,$ satisfies the {\it Figiel condition} if there exists a constant $K$ such that the function $N(\cdot)$ on some interval $(0, \e)$ satisfies the condition
\begin{equation}
\frac{N(s)}{s^2} \men K \frac{N(t)}{t^2} \quad \forall\; 0 < t \men s < \e.
\label{eq:Orlicz_Cond}
\end{equation}
\begin{remark}
\label{rem_figiel}
The modulus of smoothness of an arbitrary Banach space satisfies the Figiel condition (see\cite{Lind}, Proposition 1.e.5).)

\end{remark}
\begin{definition}
	Through $\mathfrak{M}_2$ denote the class of functions from  $\mathfrak{M}$ that satisfy the Figiel condition.
\end{definition}
It is easy to see  that an arbitrary function $f(\cdot) \in \mathfrak{M}_2$ that is not identically zero satisfies the condition
$t^2 = \operatorname{O}(f(t))$ as $t \to 0.$

Obviously, in case of a Hilbert space $H$ for the function $\psif{t} = t^2$ the assertion of Theorem \ref{th_main_hilbert} is equivalent to the equality $\Omega_{N}^{\psi}(R) = \Omega_{PS}(R).$

In the author's opinion, the condition $\psi \in \mathfrak{M}_2$ is rather natural. Moreover, it would be interesting to prove or contradict the following hypothesis.
\begin{conjecture}
	If in a uniformly convex and uniformly smooth Banach space $X$ a set
	$A \subset X$ belongs to the class $\Omega_{PS}(R)$  
	and to the class $\Omega_{N}^{\psi}(r)$ for some function $\psi \in \mathfrak{M} \setminus \mathfrak{M}_2$ and constants $R>0, r>0,$ then it belongs to the class $\Omega_{N}^{\psi_1}(r),$ where $\psi_1(t) = ct^2$ for some $c \bol 0.$
\end{conjecture}

The following two theorems provide an ample description of the connection of classes  $\Omega_N^{\psi}(R), \Omega_{PS}(R).$ 
\begin{fact}\label{Th_prox_to_npsi}
In a uniformly convex and uniformly smooth Banach space $X$  
	the following statements are equivalent for the function   $\psi \in \mathfrak{M}$:
	\begin{enumerate}
		\item  there exists  $k_1 > 0$ such that  $\Omega_{PS}(R)\subset \Omega_N^{k_1\psi}(R)$ for any $R>0;$
		\item $\rho_X(\tau) = \operatorname{O}(\psi(\tau))$ as $\tau \to 0$.
	\end{enumerate}
\end{fact}
\begin{fact}\label{Th_npsi_to_prox}
In a uniformly convex and uniformly smooth Banach space $X$  
the following statements are equivalent for the function   $\psi \in \mathfrak{M}$:
	\begin{enumerate}
		\item  there exists $k_2 > 0$ such that  $\Omega_N^{k_2\psi}(R) \subset \Omega_{PS}(R)$ for any $R>0;$
		\item $\psi(\e) = \operatorname{O}(\mcox{\e})$ as $\e \to 0$.
	\end{enumerate}
\end{fact}

Clearly, $A \in \Omega_{PS}(1)$ iff $RA \in \Omega_{PS}(R)$ for some $R > 0.$ 
A similar property holds  for the class of $\psi$-hypomonotonous sets. Therefore, Theorems \ref{Th_prox_to_npsi}, \ref{Th_npsi_to_prox} can be proved assuming that $R$ is a fixed positive number.
Using Theorems \ref{Th_npsi_to_prox}, \ref{Th_prox_to_npsi} we will show that the following statement, which gives a complete answer to the question about the possibility of the inclusion
$\Omega_N^{k_1\psi}(R)\subset\Omega_{PS}(R) \subset \Omega_N^{k_2\psi}(R)$, holds.

\begin{fact}\label{mainTh_3}
	If in a uniformly smooth and uniformly convex Banach space $X$ for some function $\psi \in \mathfrak{M}$  there exist $k_1>0$, $k_2>0$ such that the inclusions
	$\Omega_N^{k_1\psi}(R)\subset\Omega_{PS}(R) \subset \Omega_N^{k_2\psi}(R)$ hold, then $\mcox{\e}\asymp\mglx{\e}\asymp\e^2$ as $\e \to 0,$ and, therefore,
	the space $X$ is isomorphic to a Hilbert space.
\end{fact}

\begin{conjecture}
The equality $\Omega_{PS}(R)=\Omega_N^{\psi}(R)$ holds only in a Hilbert space provided that $\psi(t) = t^2.$
\end{conjecture}

\section{Preliminaries}
Geometrical properties of class $\Omega_{PS}(R)$ hide in its definition. 
To clarify the geometrical properties of $\Omega_{PS}(R)$, which are very useful in this paper,
we introduce two equivalent (in certain spaces!) definitions  of class  $\Omega_{PS}(R).$
 \begin{definition}\label{def_O_P}
 	We say that a closed set $A$ in a Banach space $X$  satisfies the
 	{\it P-supporting condition of weak convexity with constant} 
 	$R > 0$ if the facts that $u \in U(R, A)$ and $x \in P_A(u)$ 
 	imply the inequality
 	$$\rho\!\left(x + \frac{R}{\norm{u-x}}(u-x), A\right) \bol R.$$
 	The set of all such sets we denote by $\Omega_P(R).$
 \end{definition}
 \begin{definition}\label{def_O_N}
 	We say that a closed set $A$ in a Banach space $X$  satisfies the
 	{\it N-supporting condition of weak convexity with constant} 
 	$R > 0$ if the facts that $x \in A$, $ p \in N(x, A) \cap \SSS$ and
 	$u \in \SS$ such that $p \in J_1(u)$  
 	imply the inequality
 	$\rho(x + Ru, A) \bol R.$
 	The set of all such sets we denote by $\Omega_N(R).$
 \end{definition}
   The equality $\Omega_{PS}(R) = \Omega_P (R)$ was proved in \cite{IvBal} in case of
 uniformly convex and uniformly smooth Banach space. The equality $\Omega_{P}(R) = \Omega_N(R)$ was proved in \cite{IvIvSupp2011} in case of uniformly convex Banach space.  So we can formulate the following theorem.
\begin{nfact} \label{oporuslTh3.2}
Let $X$ be  a uniformly smooth and uniformly convex Banach space.
Then
$$
	\Omega_{PS}(R) =\Omega_{P}(R) = \Omega_{N}(R).
$$
\end{nfact}
 Further we will use class $\Omega_N(R)$ instead 
 of $\Omega_{PS}(R),$ as mentioned above these classes coincide under conditions of theorems \ref{Th_prox_to_npsi}, \ref{Th_npsi_to_prox}, \ref{mainTh_3}.
 
\begin{fact} \label{Th o sechenii}
Let $X$ be a reflexive Banach space, $R>0$, $A\in\Omega_{N}(R)$, 
$a_{0}\in A$, $p \in N(a_{0},A) \cap\SSS$, $\delta>0$, 
$\e  = \frac{2R}{\delta} \mglx{\frac{\delta}{R}}$. Then
$$
\bra p,a-a_{0}\ket \le \e\norm{a-a_{0}} \quad \forall a\in A\cap\B_{\delta}(a_0).
$$
\end{fact}
\begin{prf}
	Fix a  vector $a \in A \cap \B_\delta(a_0).$
	By reflexivity of $X,$ there exists a vector $u \in \SS$ such that 
	$p \in J_1(u).$ Since $A \in \Omega_N(R),$ we have that $\rho(a_0 + Ru, A) \bol R.$
	Therefore, $\norm{a_0 + Ru - a} \bol R.$ Denote $w = \frac{a -a_0}{R}.$ Then $\norm{u - w} \bol 1.$
	By the definition of the modulus of smoothness we get that $2\mglx{\norm{w}} \bol \norm{u+w} + \norm{u-w} -2.$
	Consequently,
	$$
	2\mglx{\norm{w}} \bol \norm{u+w} -1 \bol \bra p, u+w\ket -1  = \bra p, w\ket.
	$$
	And now we have 
	$$
		\bra p, a - a_0\ket = R\bra p, w\ket \men 2 R\mglx{\norm{w}} = 2R \mglx{\frac{\norm{a - a_0}}{R}}.
	$$
	Since  function $\tau \to \frac{\mglx{\tau}}{\tau}$ is increasing and $\norm{a - a_0} \men \delta,$
	we have $\mglx{\frac{\norm{a-a_0}}{R}} \men \frac{\norm{a - a_0}}{\delta} \mglx{\frac{\delta}{R}}.$
	Consequently, $\bra p,a-a_{0}\ket \le \e\norm{a-a_{0}}.$
\end{prf}

We need to formulate some properties of the moduli of smoothness and convexity.
 
It is well known (\cite{DiestelEng} Ch.3, \S 4, Lemma 1.), 
that in a Banach space the following inequalities hold.
\begin{equation}
\label{mgl2_2}
2 \men \limsup\limits_{\tau \to 0} \frac{\mglx{2\tau}}{\mglx{\tau}} \men 4.
\end{equation}
The following two lemmas are technical and obvious.
The first one can find in \cite{Supp_mod_arxiv}. 
\begin{lemma}\label{UVO lemma ozen}
Let $x,y \in X,\; x\neq 0,\; p \in J_1(x).$ Then
\begin{equation}
    \norm{x+y} \men \norm{x} + \bra p, y \ket + 2\norm{x} \mglx{{\frac{\norm{y}}{\norm{x}}}}.
\end{equation}
\end{lemma}
\begin{lemma}\label{lemma_ozenka_otkl_mcox}
Let $X$ be a uniformly convex Banach space,  $x_0 \in \partial\B_R(o)$ and  
	$p_0 \in \SSS$ be a functional dual to the vector $-x_0.$ 
	Then for any vector  $z \in \vn\B_R(o)$ we have
	\begin{equation}\label{lemma_ozenka_otkl_mcox2}
		\bra p_0, z- x_0\ket > 2 R \mcox{\frac{\norm{z-x_0}}{R}}.
	\end{equation}
\end{lemma}

The proof of  the next proposition one could find in work \cite{IvIvSupp2011}.
\begin{nfact} \label{oporuslLm3.1}
Let  $A\subset X$, $x_{1}\in X\setminus A$, 
$x_{0}\in P_{A}(x_{1})$, and the norm is Frechet differentiable at point  $x_{1}-x_{0}$. Then
$J_1(x_{1}-x_{0})\subset N(x_{0},A)$.
\end{nfact}
We will denote by $\B_1^p(o)$ the intersection of the unit ball 
and the hyperplane 
$$H_p = \left\{ x \in X | \;\bra p, x\ket = 0 \right\},$$
i.e. the set $\BB \cap H_p.$ Let $\B_1^p(a) =  \B_1^p(o) + oa.$ 
\begin{lemma}\label{lemma_o_proek_hordi}
	Let $x \in \SS$, $p \in J_1(ox),$  $z \in \B_1^p(x).$ 
	Let a point  $y$  such that segment $yz$ is parallel to  $ox$
	and intersects the unit sphere at unique point $y.$
	Then $2\norm{zx} \bol \norm{xy}.$
\end{lemma} 
\begin{prf}
	By the triangle inequality, it suffices to show that $\norm{zx} \bol \norm{zy}.$
	By   $d$ denote a point on unit sphere such that vector $od$ is collinear to the vector $xz.$
	Let line $l$ be parallel to $ox$ and $d \in l.$
	By constructions, we have that  points $x,y,z,o,d$ and line	$l$ 
	lies  at the same plane  -- linear span of the vectors $ox$ and $xz$.
	So lines $l$ and $xz$ intersect, by $c$ we denote their intersection point. 
	Note that  $odcx$ is a parallelogram and $\norm{dc} = 1$;
	segment $dx$ belongs to the unit ball and does not intersect the interior of the segment $zy.$
	Let $y' = zy \cap dx.$
 	By similarity, we have
	$$
		\norm{zy} \men \norm{z y'} = \frac{\norm{x z}}{\norm{x  c}} \norm{d  c} = \norm{x  z}. 
	$$
\end{prf}

	It is worth noticing, that  in the conditions of   Lemma \ref{lemma_o_proek_hordi} we have 
	that $z$ is a projection along vector $ox$ of the point $y$
	 on some supporting hyperplane to the unit ball at $x$. 
	Moreover, $z$  belongs to the metrics projection of the point $y$ on this hyperplane. 
	That is to say,   Lemma \ref{lemma_o_proek_hordi} show us that 
	if one projects along vector $ox$ segment $xy$ on the   hyperplane, which is supporting to the unit ball	at point $x,$ then the length of the segment decrease no more than a factor of 2
\begin{figure}[h]%
\center{
\psfrag{x}[1]{$\hspace{0.3em}x$}
\psfrag{o}[1]{$o$}
\psfrag{y}[1]{\hspace{0.3em}$y$}
\psfrag{z}[1]{\hspace{0.5em}$z$}
\psfrag{B}[1]{\raisebox{-3.5ex}{$\B_1^p(x)$}}
\psfrag{l3}[1]{\hspace{-0.8em}$l_3$}
\psfrag{l1}[1]{$l_2$}
\psfrag{l}[1]{\raisebox{2.3ex}{$l$}}
\includegraphics[scale=0.5]{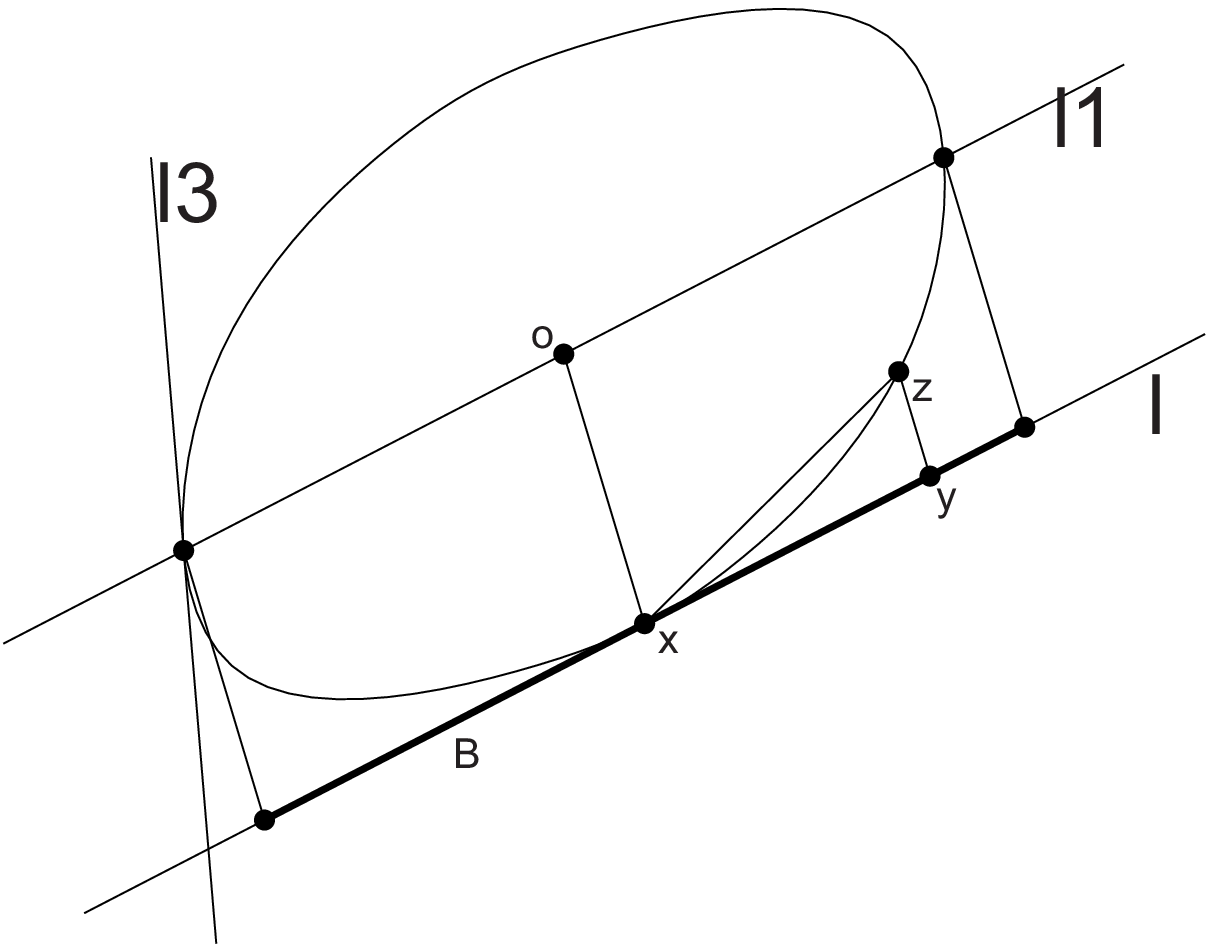}
}
\caption{ Under the conditions of Lemma \ref{lemma_o_proek_hordi} we have $2\norm{zx} \bol \norm{xy}.$}
\end{figure}

In this paper we will use the notion of the moduli of supporting convexity and smoothness, which 
were introduced and studied by the author in paper \cite{Supp_mod_arxiv}. This technique is very convenient in questions concerning
the geometry of  unit ball and its supporting hyperplanes. 
Let us introduce the following definitions.

We say that $y$ is {\it quasiorthogonal}  to vector $x \in X\setminus \{o\}$  
and write  $y\urcorner x$ 
if there exists a functional $p \in J_1(x)$  such that  $\bra p, y \ket = 0.$
Note that  the following conditions are equivalent: \\ \noindent
 -- $y$  is quasiorthogonal to  $x$ \\ \noindent 
 --  for any $\la \in \R$ the vector $x+ \la y$ 
 lies in the supporting hyperplane to the  ball  $\B_{\norm{x}}(o)$  at $x;$ \\ \noindent
 -- for any $\la \in \R$ the following inequality holds  $\norm{x + \la y} \bol \norm{x};$\\ \noindent
 --  $x$ is orthogonal to  $y$ in  Birkhoff--James sense (\cite{DiestelEng}, Ch. 2, \S 1).
 
Let $x,y \in \SS$ be such that  $y\urcorner x.$ 
By definition, put  
$$\la_X(x,y,r) = \min{\{\la \in \R: \, \norm{x+ry - \la x} = 1\}}$$
for any  $r \in [0,1].$ 

Define the {\it modulus of supporting convexity} as
$$
\lambda^{-}_{X}(r) = \inf \{\la_X(x,y,r):\, \|x\| = \|y\|=1,\,
 y \urcorner x \}.
$$

Define the {\it modulus of supporting smoothness} as
$$
\lambda^{+}_{X}(r) = \sup \{\la_X(x,y,r):\, \|x\| = \|y\|=1,\, y \urcorner x \}.
$$

In paper \cite{Supp_mod_arxiv} the following inequalities are proved:
\begin{equation}\label{mglx_suppmglx}
\mglx{\frac{r}{2}} \men \lapx{r} \men \mglx{2r}.
\end{equation}
\begin{equation}\label{suppconvprop1}
\lamx{r} \men  \mcox{2r} \men 1 - \sqrt{1 - r^2}.
\end{equation}
\begin{equation}\label{suppconvprop2}
0 \men \lamx{r} \men  \lapx{r} \men r.
\end{equation}

\section{Proof of Theorem \ref{Th_prox_to_npsi}}
\noindent {\bf Auxiliary results}\\ \noindent

For the set $A \subset X$ we introduce the following function
\begin{equation*}
\Gamma(A, \e, X) = -\inf \bra p_1 - p_2, x_1 - x_2\ket,
\end{equation*}
where the infimum is considered for all $x_1, x_2$ in $X$ and  $p_1, p_2$ in $X^*$ 
that satisfy the following condition
\begin{equation}
x_i\in A,\; p_i \in N(A, x_i) \cap \SSS,\; i=1,2,
\label{eq:hypo_ball_cond}
\end{equation}
and such that $\norm{x_1-x_2} = \e.$

The function  $-\Gamma(A, \e, X)$ is the maximum value that can be on the right hand side in the hypomonotonity inequality
\reff{cone_to_prox_hypo1} meaning that if
$A \in \Omega_N^{\psi}(R),$ then
\begin{equation}
\label{eq:thmain1_15.46}
- \Gamma(A, \e, X) \bol -R \psif{\frac{\e}{R}} .
\end{equation}

\begin{lemma}\label{hypo_ball_1}
Let $X$ be a Banach space with a Frechet differentiable norm
Let $A = X \setminus \vn \B_1(o).$ Then
	\begin{equation}
	\Gamma(A, \e, X) \asymp \lapx{\e} \asymp \mglx{\e}, \quad \e \to 0.
	\label{eq:lemma_hypo_ball}
	\end{equation}
	Moreover, if $\e \in [0, \frac{1}{2}]$ then the inequality
	\begin{equation}
	\Gamma(A, \e, X) \bol \mglx{\frac{\e}{4}}.
	\label{eq:lemma_hypo_ball_1-1}
	\end{equation}
	holds.
\end{lemma}
\begin{prf}
1) The equivalence of $\lapx{\e} \asymp \mglx{\e}$ as $\e \to 0$  is the consequence of \reff{mglx_suppmglx}.\\ \noindent
2) Fix any $\e \in (0, \frac{1}{2}].$ 
As $\partial A = \partial\B_1(o),$ 
it follows that for any $x \in \partial A$  the equality
\begin{equation}
	N(A, x) \cap \SSS = - J_1(x)
\label{eq:lemma_hypo_ball_2}
\end{equation}
holds.
Fix any $x_1,x_2, p_1, p_2,$ that satisfy the inclusions \reff{eq:hypo_ball_cond} and such that 
$\norm{x_1 -x_2} \men \e \men \frac{1}{2}.$
Let $\ell_1$ be a line in the plane $x_1ox_2,$ supporting to the ball $\BB$ at the point $x_1$
(its uniqueness follows from the differentiability of the norm). 
Let $z$ be the projection of point $x_2$ on this line.
Then $- \bra p_1, x_1 -x_2 \ket = \norm{z-x_2}.$
From the triangle inequality it follows that
$
\norm{z-x_1} \men \norm{z-x_2} + \norm{x_2- x_1} \men 2 \norm{x_2 - x_1}.	
$
Therefore, $$\norm{z-x_2} \men \lapx{\norm{x_1 - z}} \men \lapx{2\norm{x_1 - x_2}}.$$
From the last inequality in this series we obtain that
$
\bra p_1, x_1 - x_2 \ket \bol - \lapx{2\norm{x_1 - x_2}}.
$

Similarly, $\bra p_2, x_2 - x_1 \ket \bol - \lapx{2\norm{x_1 - x_2}}.$
Summing the last two inequalities, we obtain that
\begin{equation*}
	\bra p_1 - p_2 , x_1 - x_2 \ket \bol - 2\lapx{2\norm{x_1 - x_2}},
\end{equation*}
and, so, for any $\e \in (0, \frac{1}{2}]$ we have
\begin{equation}
	\Gamma(A, \e, X) \men 2 \lapx{2\e}.
	\label{eq:hypo_ball_3}
\end{equation}

The equality \reff{eq:lemma_hypo_ball_2} implies that for any $x_1, x_2, p_1, p_2,$ 
that satisfy the condition  \reff{eq:hypo_ball_cond} the inequality
\begin{equation}
	\bra p_1 - p_2, x_1 -x_2\ket \men \bra p_1, x_1 - x_2\ket
	\label{eq:hypo_ball_4}
\end{equation}
holds.
According to the definition of the convexity modulus, for any $\gamma > 0$ 
there exist  $x_1, x_2 \in \partial A;$ $p_1 \in N(x_1, A_1);$
$z\in H_{p_1} = \left\{ x \in X |\; \bra p_1, x \ket= 1\right\}$ 
such that
\begin{enumerate}
	\item $z$ is the projection of $x_2$ on the hyperplane $H_{p_1};$
	\item $\norm{x_1 - z} = \frac{\e}{2};$ 
	\item $\norm{z-x_2} \bol \lapx{\norm{x_1 - z}} - \gamma = \lapx{\frac{\e}{2}} - \gamma.$
\end{enumerate}
This implies that
\begin{equation}
	- \bra p_1, x_1 - x_2\ket = \norm{z-x_2} \bol \lapx{\frac{\e}{2}} - \gamma.
	\label{eq:hypo_ball_4.5}
\end{equation}
On a unit sphere in the plane $x_1ox_2$ we fix a point $x_{\e}$ such that
$\norm{x_1 - x_{\e}} = \e$ and in the plane $x_1ox_2$ the triangle $x_1x_2x_{\e}$ lies on one side of line $ox_1.$
Lemma \ref{lemma_o_proek_hordi} implies that $\e = 2 \norm{x_1-z} \bol \norm{x_1-x_2}.$
Using the inequality \reff{eq:hypo_ball_4}, we obtain 
$$\Gamma(A, \e, X) \bol - \bra p_1, x_1 - x_{\e}\ket \bol - \bra p_1, x_1 - x_2\ket.$$ 
This and inequality \reff{eq:hypo_ball_4.5} imply that
\begin{equation*}
	\Gamma(A, \e, X) \bol  \lapx{\frac{\e}{2}} - \gamma.
\end{equation*} 
Passing to the limit in the last inequality as $\gamma \to 0,$ we obtain that
\begin{equation}
	\Gamma(A, \e, X) \bol \lapx{\frac{\e}{2}}.
	\label{eq:hypo_ball_5}
\end{equation}
The inequalities \reff{eq:hypo_ball_3} and \reff{eq:hypo_ball_5} imply that
\begin{equation}
	2 \lapx{2\e} \bol \Gamma(A, \e, X) \bol \lapx{\frac{\e}{2}}
	\label{eq:hypo_ball_6}
\end{equation}
Using this and \reff{mglx_suppmglx} we obtain inequality \reff{eq:lemma_hypo_ball_1-1}.
\end{prf}

\noindent {\bf Proof $2 \Rightarrow 1$ in Theorem \ref{Th_prox_to_npsi}}\\ \noindent 
Fix a set $A \in \Omega_N(R).$
Theorem \ref{Th o sechenii} implies that for any $x_1, x_2 \in A$ и $p_i \in N(x_i, A), i =1,2$ the inequality
\begin{equation}
\bra p_1 - p_2 , x_1 - x_2\ket \bol  - 4 R \mglx{\frac{\norm{x_1 - x_2}}{R}}
\label{cone_to_prox:2}
\end{equation}
holds.
Condition 2 implies that for some $\e >  0$ exists a positive constant $k$ such that
for any $\tau \in [0, \e]$ the inequality
\begin{equation}
	\mglx{\tau} \men k \psif{\tau}
	\label{cone_to_prox:1106.1}
\end{equation}
holds.
The inequalities \reff{cone_to_prox:1106.1}, \reff{cone_to_prox:2} imply that
$A \in \Omega_N^{k_1\psi}(R),$ where $k_1 = 4k.$
\bbox

\noindent {\bf Proof $1 \Rightarrow 2$ in Theorem \ref{Th_prox_to_npsi}}\\ \noindent 
Let the constant $k_1$ be such that $\Omega_{N}(1)\subset \Omega_N^{k_1\psi}(1).$
Let us show that assertion $2$ holds. 
Consider the set $A_1 = \{x \in X | \; \norm{x} \bol 1\}.$
Clearly, $A_1 \in \Omega_P(1),$ and by Theorem \ref{oporuslTh3.2} the inclusion $A_1 \in \Omega_N(1)$ holds.
From the definition of $\Gamma(A, \e, X),$ inequalities \reff{eq:lemma_hypo_ball_1-1}, \reff{eq:thmain1_15.46} and the relations \reff{mgl2_2} 
we obtain for some constant $\gamma > 0$ the following series of inequalities inequalities
\begin{equation*}
	k_1 \psif{{\e}} \bol \Gamma(A, \e, X) \bol \mglx{\frac{\e}{4}} \bol \frac{1}{17} \mglx{\e}, 
\end{equation*}
where $\e \in [0, \gamma].$

\bbox 
Particularly, we have shown that the following statement holds. 
\begin{corollary}\label{corollary_o_perenorm}
In a uniformly smooth and uniformly convex Banach space $X$ the inclusion
	$\left(X \setminus \vn \BB\right) \in \Omega_N^{\frac{1}{17}\mglx{\cdot}}(1)$
	holds.
\end{corollary}
\begin{remark}
	The implication $2 \Rightarrow 1$ in Theorem \ref{Th_prox_to_npsi} is an elementary corollary of the smoothness modulus definition 
	(see Theorem \ref{Th o sechenii}). But to prove the reverse implication we had to show the equivalence of smoothness and convexity moduli at zero. 
\end{remark}

\section{\bf Proof 2 $\Rightarrow$ 1 in Theorem \ref{Th_npsi_to_prox}}
\noindent {\bf Auxiliary results} \\ \noindent
Define the function $\omega_X:  [0,+\infty) \to  [0,+\infty)$ in the following way
$$
	\omega_X(\tau) = \frac{\mglx{\tau}}{\tau}.
$$
As  $\omega_X(\cdot)$ is an increasing monotonous function, the inverse function $\omxi{\cdot}$ is also increasing.
Note that $\omega_X(1) = \mglx{1} \bol \mgl{H}{1} = \cmgl,$ 
and, therefore, for any $t \in [0, 1]$ the inequality
\begin{equation}
	\omxi{(\cmgl) t} \men 1,
\label{ineq_omega}
\end{equation} 
holds.

The next Lemma is a generalization of Lemma 1.9.1 from G.E. Ivanov's monograph \cite{IvMonEng},
in which an analogous result was obtained for Hilbert spaces.

For a set $A \subset X$ denote by $T(A)$ the set of a point $x$ such that $P_A(x)$ is single-point.
It is well known, that in a uniformly convex Banach space  $X$ for any closed set $A$ the set $T(A)$ is dense 
in $X$ (See \cite{Stechkin_TH}).
\begin{lemma}\label{lemma_cone_to_prox1}
	Let $X$ be a uniformly smooth Banach space, 
	$A$ be a closed set such that the set $T(A)$ is dense in $X$.
	Let there be given $z_0 \notin A,\, z_1 \in A$ and $\e \in (0,1).$
	Then there exist $\la \in [0,1],$ a point  $y \in \partial A$ and $p \in N(y, \, A), \norm{p} = 1$ such that
	\begin{equation}\label{prox eq:1}
	\norm{(1-\la)z_0 + \la z_1 - y} < \e \norm{z_1-z_0},
	\end{equation}
	\begin{equation}\label{prox eq:2}
	\bra p, z_1 - z_0 \ket <  \e \norm{z_1-z_0}.
	\end{equation}
\end{lemma}
\begin{prf}
	Note
	$
	\dd = \frac{1}{2} \min\{ \e \norm{z_1 - z_0},\, \rho(z_0,\,A)\}.
	$
 As the set $A$ is closed and $z_0 \notin A,$ we have that $\rho(z_0,\,A) >0,$ and, therefore $\dd >0.$
	Since $2\dd \men \rho(z_0,\,A)$, we obtain that $\vn \B_{2\dd}(z_0) \subset X \setminus A.$
	
	Define the number $\la$ as follows
	\begin{equation*}
	\la  = \inf \left\{t > 0: \B_\dd\left( z_0 + t(z_1 - z_0)\right) \bigcap A \neq \emptyset \right\}.
	\end{equation*}

	As $z_1 \in A,$ we have that $\B_\dd(z_1) \cap A \neq \emptyset,$ and, therefore, $\la \men 1.$
	Note that for  $|t| < \frac{\dd}{\norm{z_0 -z_1}}$ the inclusions
	\begin{equation*}
	\B_\dd\left( z_0 + t(z_1 - z_0)\right) \subset \vn \B_{2\dd}(z_0) \subset X \setminus A,
	\end{equation*}
	hold and, therefore,  $\B_{\dd}\left( z_0 + t(z_1 - z_0)\right) \cap A = \emptyset$ as
	$|t| < \frac{\dd}{\norm{z_0 -z_1}}.$
This and the definition of number  $\la$ imply that
	\begin{equation}
	\la \bol \frac{\dd}{\norm{z_0 -z_1}},
	\label{prox eq:3}
	\end{equation}
	\begin{equation}
	\B_\dd\left( z_0 + t(z_1 - z_0)\right) \cap A = \emptyset,\;\; \forall t \in\left( -\frac{\dd}{\norm{z_0 -z_1}}, \la  \right).  
	\label{prox eq:4}
	\end{equation}
	
	Define the vector $z_{\la} = (1-\la)z_0 + \la z_1.$ 
	Since $X$ is a uniformly smooth Banach space, there exists $\xi > 0$ such that
	\begin{equation}
	(\sqrt{2} + 1)\frac{\xi}{\mgli{X}{\xi}} < \e
	\label{prox eq:5}
	\end{equation} 
	Define the number  $\beta = \frac{\sqrt{2}-1}{2}\dd \min\{\xi, \frac{1}{2}\}.$
	
	As the set $T(A)$ is dense in $X$, there exists $y_0 \in T(A)$ such that
the inequality 
	$\norm{z_{\la} - y_0} < \beta$ holds.
	Let  $\{y\} = P_A(y_0).$

	Obviously,
	\begin{equation*}
		\left|\rho(y_0, A) - \rho(z_{\la}, A)\right| \men \norm{y_0 - z_{\la}} < \beta.
	\end{equation*}
	By the definition of  $\la$ we have $\rho(z_{\la}, A) = \delta.$
	Therefore,
	\begin{equation}
		\dd - \beta < \norm{y_0 - y} =  \rho(y_0, A)< \dd +\beta.
		\label{prox eq:6}
	\end{equation} 
	From this and the inequality $\norm{y_0 -z_\la} < \beta,$ $\beta \men \frac{\dd}{2},$ $\dd \men \frac{\e\norm{z_1 - z_0}}{2}$
	we get that
	\begin{equation*}
		\norm{z_\la - y} \men \norm{y_0 - y} + \norm{y_0 - z_\la} < \dd + 2\beta \men 2\dd \men \e\norm{z_1 -z_0}.
	\end{equation*}
	Thus, inequality \reff{prox eq:1} holds.
	
Let the functional  $p \in \partial \B^*_1(o)$ be dual to the vector  $y_0 - y.$
	Lemma \ref{oporuslLm3.1} implies that $p \in N(y,\, A).$
	
	If $\bra p, z_1 -z_0\ket \men 0,$ then inequality \reff{prox eq:2} holds.
	So, let us consider the case $\bra p, z_1 -z_0\ket > 0.$
	Let $\gamma = \frac{\bra p, z_1 -z_0\ket}{\norm{z_1 -z_0}} \in (0, 1].$
	We have to prove that $\gamma < \e.$
Define the number
	\begin{equation} \label{mu:1}
		\mu = \omxi{(\cmgl) \gamma} \frac{\norm{y-y_0}}{\norm{z_1 -z_0}}.
	\end{equation}
Applying to the last expression inequalities \reff{ineq_omega},  \reff{prox eq:3} and \reff{prox eq:6}, we obtain that
\begin{equation*}
0 < \mu \men \frac{\norm{y-y_0}}{\norm{z_1 -z_0}} \men \frac{\dd + \beta}{\norm{z_1 -z_0}} < \frac{2\dd}{\norm{z_1 -z_0}} 
\men \frac{\dd}{\norm{z_1 -z_0}} + \la.
\end{equation*}
Thus,  $-\frac{\dd}{\norm{z_1 -z_0}} < \la - \mu < \la.$
This and equality \reff{prox eq:4} imply that
\begin{equation*}
\B_\dd (z_0 + (\la - \mu)(z_1- z_0)) \cap A = \emptyset.
\end{equation*}
As $y \in A,$ we have that
$y \notin \B_\dd (z_0 + (\la - \mu)(z_1- z_0)) = \B_\dd ( z_\la - \mu(z_1- z_0)),$
and so  $\norm{z_\la -y -  \mu(z_1- z_0)} > \dd. $
This and the inequality $\norm{y_0 - z_\la} < \beta$ imply that
$\norm{y_0 - y  - \mu(z_1 - z_0)} > \dd - \beta.$
Applying Lemma \ref{UVO lemma ozen}, we obtain that
\begin{equation*}
\frac{\dd - \beta}{\norm{y_0 - y}} < 1 - \mu\frac{ \bra p , z_1- z_0\ket}{\norm{y_0-y}}  + 2 \mglx{\mu \frac{\norm{z_1 -z_0}}{\norm{y_0-y}}}.
\end{equation*}
Thus, considering inequalities  \reff{mu:1}, \reff{prox eq:6} we get that
\begin{equation*}
\frac{\dd - \beta}{\dd + \beta} < 1 - \gamma \omxi{(\cmgl) \gamma}+ 2 \mglx{\omxi{(\cmgl) \gamma}}.
\end{equation*}
Note that for any  $k > 0$ the equalities
$$\gamma \omxi{k \gamma} = \frac{\omxi{k\gamma} \cdot \omega_X\left( \omxi{k\gamma}\right)}{k} = \frac{\mglx{\omxi{k\gamma}}}{k}$$ 
hold,
and, therefore, we get that
\begin{equation*}
\frac{\dd - \beta}{\dd + \beta} < 1 - \left( \frac{1}{\cmgl} - 2\right)  \mglx{\omxi{(\cmgl) \gamma}}.
\end{equation*}
Modifying  the last inequality and using the definition of 
$\beta,$ we obtain that
\begin{equation*}
\mglx{\omxi{(\cmgl) \gamma}} < \frac{2(\sqrt{2}-1)}{(3-2\sqrt{2})} \frac{\beta}{(\delta+\beta)}= \frac{2}{(\cmgl)} \frac{\beta}{(\delta + \beta)} \men \xi,
\end{equation*}
and thus, considering inequality  \reff{prox eq:5}, we get that
\begin{equation*}
\gamma < \frac{1}{\cmgl}\omega_X\left(\mgli{X}{\xi}\right)  =
\frac{1}{\cmgl} \frac{\mglx{\mgli{X}{\xi}}}{\mgli{X}{\xi}}
=  (\sqrt{2} + 1)\frac{\xi}{\mgli{X}{\xi}} < \e.
\end{equation*}
\end{prf}

\noindent {\bf Proof of the main assertion}\\ \noindent
	As the functions $\mcox{\cdot}$ and $\psif{\cdot}$ are monotonous, condition 2 of theorem \ref{Th_npsi_to_prox} implies that there exists a number $k > 0$ such that
	\begin{equation}
	\label{cone_to_prox_rem_ofunc}
	R \cdot k \psif{\frac{t}{R}}\men 	2R \cdot\mcox{\frac{t}{R}} \quad \frl t \in [0, 2R].
	\end{equation}
Let us prove that $\Omega_N^{k\psi}(R) \subset \Omega_N(R).$

Suppose the contrary.
Then there exists a closed set $A$ such that $A \in \Omega_N^{k\psi}(R)$ and $A \not \in \Omega_N(R).$  
Then there exist a vector $x_0 \in A,$  a functional $p_0 \in N(x_0, A) \cap \partial \B_1^{*}(o)$ and the dual to it unit vector $u_0,$ and also the vector  $z_1 \in A$
such that the inequality
\begin{equation}
\norm{x_0 + R u_0 - z_1} < R
\label{eq:lemma_ozenka_otkl_23503009}
\end{equation}
holds.

According to Lemma \ref{lemma_ozenka_otkl_mcox}, we obtain that
\begin{equation}
\bra p_0 , z_1 - x_0 \ket >2R \mcox{\frac{\norm{z_1 - x_0}}{R}} >  0.
\label{cone_to_prox_eq1}
\end{equation}
Therefore, due to the inclusion $p_0 \in N(x_0, A)$ and according to Definition \ref{def_nconus}, 
there exists a number $\mu \in (0,1)$ such that the vector
\begin{equation}
\label{lemma_ozenka_otkl_eq2}
	z_0 = x_0 + \mu (z_1 - x_0)
\end{equation}
satisfies the condition $z_0 \not \in A.$

Let $L$ be the Lipschitz constant of function $\psif{\cdot}.$ 
Denote
 \begin{equation}
\label{lemma_ozenka_otkl_eq3}
	\e = 
	\frac{\mu}{(6+2kL)R}
	\left(\bra p_0 , z_1 - x_0 \ket - 2R \mcox{\frac{\norm{z_1 - x_0}}{R}} \right).	
\end{equation}
Inequality (\ref{cone_to_prox_eq1}) implies that $\e > 0.$

According to Lemma \ref{lemma_cone_to_prox1} there exists a number $\la \in [0,1],$ a vector $y \in \partial A$ 
and a unit functional $p \in N(y, A)$ such that
\begin{equation}\label{lemma_ozenka_otkl_eq4}
\norm{z_\la - y} < \e \norm{z_1-z_0},
\end{equation}
\begin{equation}\label{lemma_ozenka_otkl_eq5}
\bra p, z_1 - z_0 \ket <  \e \norm{z_1-z_0},
\end{equation}
where
\begin{equation}\label{lemma_ozenka_otkl_eq6}
z_\la = (1-\la)z_0 + \la z_1.
\end{equation}
The definition of vector $z_0$ (equality (\ref{lemma_ozenka_otkl_eq2})) implies that
$$
	z_1 - z_0 = (1 - \mu)(z_1-x_0), \quad z_0 - x_0 = \mu(z_1 - x_0).
$$
From this and equality (\ref{lemma_ozenka_otkl_eq6}) we get that
\begin{equation}
	z_{\la} -x_0 = z_\la -z_0 + z_0 -x_0 
	=\la(z_1 -z_0)+z_0-x_0 = (\la(1-\mu) + \mu)(z_1 - x_0). \label{lemma_ozenka_otkl_eq7}
\end{equation}
Therefore
$$
	\frac{z_\la - x_0}{\norm{z_\la - x_0}} = \frac{z_1 - x_0}{\norm{z_1 - x_0}} = \frac{z_1 - z_0}{\norm{z_1 - z_0}}.
$$
This and inequality (\ref{lemma_ozenka_otkl_eq5}) imply that
$$
	\bra p, z_\la - x_0\ket < \e \norm{z_\la - x_0} \men \e \norm{z_1 - x_0}. 
$$
Thus, according to inequality \reff{lemma_ozenka_otkl_eq4}, 
we have that $\bra p, y - x_0 \ket \men 2\e \norm{z_1 - x_0}.$
Inequality \reff{eq:lemma_ozenka_otkl_23503009} implies that
$\norm{z_1 - x_0} \men 2R,$ and, therefore,
\begin{equation}\label{lemma_ozenka_otkl_eq8}
\bra p, y - x_0 \ket \men 4\e R. 
\end{equation}
As $A \in \Omega_N^{k\psi}(R),$  the inequality
$$
\bra p - p_0, y - x_0 \ket \bol - Rk  \psif{\frac{\norm{y - x_0}}{R}}
$$
holds.
From this and the inequality \reff{lemma_ozenka_otkl_eq8} we get that
\begin{equation}\label{lemma_ozenka_otkl_eq9}
\bra p_0, y - x_0 \ket \men Rk\psif{\frac{\norm{y-x_0}}{R}} + 4\e R.
\end{equation}

The relations \reff{lemma_ozenka_otkl_eq7} imply that $\norm{z_\la - x_0} \men \norm{z_1 - x_0} <2 R.$
Therefore, by inequality \reff{lemma_ozenka_otkl_eq4}, we have that
$$
	\norm{y-x_0} \men \norm{z_\la - x_0} + \e \norm{z_1 - z_0} < \norm{z_\la - x_0} +  2\e R.
$$

Thus, considering the Lipschitzness of function $\psif{\cdot},$ we get that
$$
	 Rk\psif{\frac{\norm{y-x_0}}{R}} \men Rk\psif{\frac{\norm{z_\la-x_0}}{R}} + 2\e L kR.
$$
Substituting the last inequality in  \reff{lemma_ozenka_otkl_eq9},  we obtain that
$$
	\bra p_0, y - x_0 \ket \men    Rk\psif{\frac{\norm{z_\la-x_0}}{R}} + \e R(4 + 2kL).
$$
Thus, using inequality \reff{lemma_ozenka_otkl_eq4}, we get that
$$
	\bra p_0, z_\la - x_0 \ket <   Rk\psif{\frac{\norm{z_\la-x_0}}{R}} + \e R(6 + 2kL).
$$

Denote $t = \la (1 - \mu) + \mu.$ Then $0 < \mu \men t \men 1$ and, according to equality
(\ref{lemma_ozenka_otkl_eq7}), we get that $z_\la - x_0 = t(z_1 - x_0).$ Thus,
$$
	t \bra p_0, z_1 - x_0 \ket <   Rk\psif{t\frac{\norm{z_1-x_0}}{R}} + \e R(6 + 2kL),
$$
due to the convexity of function $\psif{\cdot}$, we obtain that
$$
	\bra p_0, z_1 - x_0 \ket <  Rk \frac{\psif{t\frac{\norm{z_1-x_0}}{R}}}{t} + \frac{\e R(6 + 2kL)}{t} 
	\men R k\psif{\frac{\norm{z_1-x_0}}{R}} + \frac{\e R(6 + 2kL)}{\mu}.
$$
Considering inequality \reff{cone_to_prox_rem_ofunc}, we get that
$$
\bra p_0, z_1 - x_0 \ket < 2R \mcox{\frac{\norm{z_1-x_0}}{R}} + \frac{\e R(6 + 2kL)}{\mu},
$$
which contradicts equality \reff{lemma_ozenka_otkl_eq3}.
\par\noindent\ensuremath{\Box}\par\noindent

\section{Proof 1 $\Rightarrow$ 2 in Theorem \ref{Th_npsi_to_prox}}
The proof consists in constructing examples of sets in the class
 $\Omega_N^{k_2\psi}(R)$  that do not belong to the class
 $\Omega_N(R)$, and consists in number of steps.
The main idea is building in the plane $L \subset X$ a set $A_1$ that belongs to the class $\Omega_N^{k\psi}(R)$ in space $L.$
Then, with the help of some technical lemmas, the set $A_1$ extends to the set $A,$ which already is in the class $\Omega_N^{K\psi}(K_1R)$ in $X,$
where the constants $k, K, K_1$ do not depend on the initially chosen plane $L$ and space $X.$ 
Using the fact that $\mcox{t} \asymp \lamx{t} = \omal{\psif{t}}$ at zero, 
we chose such vectors that their existence contradicts the  $N$-supporting condition.

\noindent {\bf Auxiliary results} \\ \noindent
Let us first describe the functions $\psi$ from the class $\mathfrak{M}$, such that  there exists a normed space whose smoothness modulus is equivalent to the function   $\psi.$   
In paper \cite{modulsmooth1} the following theorem is proved. 
\begin{nfact}\label{NTh_o_mgl}
For the function $N : [0, +\infty) \to [0, +\infty),$  that satisfies the Figiel condition and such that
$N(0) = 0$ there exists a two-dimensional space $X_2$ whose smoothness modulus $\mgl{X_2}{t}$ is equivalent to $N(t)$ at zero.
\end{nfact}
The Day-Nordlander Theorem (see \cite{DiestelEng}) implies that if for  function $\psi$ there exists a Banach space with the smoothness modulus, equivalent to $\psi$ at zero,
then $t^2 = \operatorname{O}(\psif{t}).$
However, the function $\psif{\cdot}$ from $\mathfrak{M}$ such that $t^2 = \operatorname{O}(\psif{t})$ may not satisfy the Figiel  condition (see \reff{eq:Orlicz_Cond}), but the following Lemma holds.
\begin{lemma} \label{lemma_o_perenorm}
For any function $\psi \in \mathfrak{M}$ such that $t^2 = \omal{\psif{t}}$ at zero 
there exists a function $\psi_1(\cdot): [0, +\infty) \to [0, +\infty)$ that satisfies the Figiel condition and such that
$t^2 = \omal{\psi_1(t)}$ at zero and $\psi_1(t)= \omal{\psif{t}}$ at zero.
 \end{lemma}
\begin{prf}
 Clearly, it is sufficient to examine the case $\psif{t} = \omal{t}$ at zero, otherwise 
 $\psi_1(t)= t^{\frac{3}{2}}.$

Consider the function $h(t) = \frac{t^2}{\psif{t}}$ on the interval $[0, 1].$ The function is continuous on this interval  and for some constant  
 $k > 0$ the inequality $h(t) \bol kt$ holds.

As the set $K=\{(t,y)|\ t\in[0,1],\ 0\le y\le h(t)\}$ is compact, its convex hull $\co K$ is also a compact set.

Define a continuous function $\gamma(t)=\max\{y\in\R|\ (t,y)\in\co K\}$, $t\in[0,1].$
The function $\gamma(\cdot)$ is concave on the interval $[0,1],$ $\gamma(0) = 0.$ Therefore, in some neighborhood of zero  $\gamma(\cdot)$ strictly increases.  
 
Define the function $\psi_1(t)=\frac{t^2}{\sqrt{\gamma(t)}}$. As $\frac{\psi_1(t)}{\psi(t)}=\frac{h(t)}{\sqrt{\gamma(t)}}\le\sqrt{\gamma(t)}\to 0$ and $\frac{t^2}{\psi_1(t)}=\sqrt{\gamma(t)}$ approaches zero monotonically as $t\to 0$, then the function $\psi_1(\cdot)$ satisfies the required conditions.
\end{prf}

Now let us prove the technical lemmas, with the help of which we will extend the set $A_1,$
which lies in the class  $\Omega_N^{k\psi}(R)$ on a plane, to the set $A,$ 
	which already lies in the class $\Omega_N^{K\psi}(K_1R)$ in the space, which contains the plane, 
	where the constants $k, K, K_1$ do no depend on the choice of the initial plane and the space containing it.
	
If the Banach space can be represented as the direct sum of two of its closed subspaces  $Z = Z_1 \oplus Z_2,$ 
then any vector $z \in Z$  is uniquely expressible as the sum of two vectors $z_1 \in Z_1$ and $z_2 \in Z_2$, and in this case we will write
$z=(z_1, z_2).$ Hereafter, when we speak about the representation of a Banach space as the direct sum of its subspaces we will consider that the subspaces are closed.
\begin{lemma} \label{lemma_o_proektore}
In the Banach space  $X$  fix any $x \in \SS,$ $p \in J_1(x),$ 
$y \in \Ker p \cap \SS$ (i.e. $y \prp x$). Denote $L = \lin{x,y}, l = \lin{y}.$ Then there exists a projector $P: X \to L$ such that
$P(\Ker{p}) = l$ and $\norm{P}\men 3.$ \nocite{kadecodopoln1}
\end{lemma}
\begin{prf}
If $\dim X = 2,$ then the identity transformation satisfies the required condition. 

Let 	$\dim X > 2.$
According to the paper  \cite{kadecodopoln2}, in the subspace $\Ker p$ of co-dimension 1 there exists a projector $P_1: \Ker p \to l$ with unit norm.
	Clearly,  $X$ is the direct sum of the subspaces $\Lin{x}$ and  $\Ker p.$ 
	
	Define the projector $P: X \to L$ in the following way
	\begin{equation*}
	Pz = P(z_1, z_2) = (z_1, P_1z_2), \qquad \frl z=(z_1, z_2) \in X,\ z_1\in \Lin x,\ z_2\in\Ker p.
	\end{equation*}
	Let us estimate its norm. Fix an arbitrary vector $u = (u_1,u_2) \in \SS,$ where $u_1\in \Lin{x}$, $u_2\in\Ker p.$
	As $u_2 \prp u_1,$ then $\norm{u_1} \men \norm{u_1 + u_2}=\norm{u} = 1.$
	From this and the triangle inequality we get  that
	$\norm{u_2} \men \norm{u_1+u_2} + \norm{u_1} \men 2.$ Therefore, 
	$$
	\norm{Pu} = \norm{(u_1, P_1u_2)} = \norm{u_1 + P_1 u_2} \men 
	\norm{u_1} + \norm{P_1u_2} \men \norm{u_1} + \norm{u_2} \men 3.
	$$ 
\end{prf}
\begin{lemma}\label{lemma_o_perenorm_2}
	Let $X$ be a linear vector space  with two equivalent norms $\norm{\cdot}_1$ and $\norm{\cdot}_2:$
	\begin{equation}
	\norm{x}_1 \men \norm{x}_2 \men c_1\norm{x}_1\!.  
 	\label{eq:lemma_o_perenorm_1}
	\end{equation}
Denote $X_1 = \left(X, \norm{\cdot}_1 \right),$  $X_2 = \left(X, \norm{\cdot}_2 \right).$ 
	Let   $A = X_2 \setminus \vn \BB$  and $A \in \Omega_N^{\psi}(R)$ in the space $X_2,$
	then $A \in \Omega_N^{c_1^2 \psi}\!\!\left(\frac{R}{c_1}\right)$ in the space $X_1.$
\end{lemma}
\begin{prf}
	According to Remark \ref{remark_def_N}, at every point of set $A$ the normal cone does not change, but the norms of the corresponding functionals may change. 
	Fix  arbitrary  $x_1, x_2, p_1^1, p_2^1,$ which satisfy the condition \reff{eq:hypo_ball_cond} for the space $X_1.$
	For $j=1,2$ denote $p_j^2=\frac{p_j^1}{\norm{p_j^1}_2}$. As $\norm{p_j^1}_1=1$, then by  inequalities \reff{eq:lemma_o_perenorm_1} we have 
	$\frac{1}{c_1}\men \norm{p_j^1}_2 \men 1$.
As the set $A$ is a complement of a convex set, then $\bra p_1^1, x_1 - x_2\ket \men 0.$
Therefore,
	\begin{equation*}
	\bra p_1^2, x_1 - x_2  \ket \bol \bra p_1^1, x_1 - x_2\ket \bol c_1 \bra p_1^2, x_1 - x_2  \ket.
	\end{equation*}
Similarly,
	 $
	\bra p_2^2, x_2 - x_1  \ket \bol \bra p_2^1, x_2 - x_1\ket \bol c_1 \bra p_2^2, x_2 - x_1  \ket.
	$
	  Summing the last two series of inequalities, we obtain that 
	\begin{equation}
	\bra p_1^2 - p_2^2 , x_1 - x_2  \ket \bol \bra p_1^1 - p_2^1, x_1 - x_2\ket \bol c_1 \bra p_1^2 - p_2^2, x_1 - x_2  \ket.
	\label{eq:lemma_o_perenorm_3}
	\end{equation}
	From the fact that $A \in \Omega_N^{\psi}(R)$ in the space $X_2$ and from the inequalities $\reff{eq:lemma_o_perenorm_1}$ we have that
 \begin{equation*}
		 c_1 \bra p_1^2 - p_2^2, x_1 - x_2  \ket 		\bol -c_1 R\psif{\frac{\norm{x_1-x_2}_2}{R}}\bol -c_1 R\psif{\frac{c_1\norm{x_1-x_2}_1}{R}}.
	\end{equation*}
	This and the last inequality in sequence \reff{eq:lemma_o_perenorm_3} imply that
	$$\bra p_1^1 - p_2^1 , x_1 - x_2  \ket \bol -c_1^2 R\psif{\frac{c_1\norm{x_1-x_2}_1}{R}}\!\!.$$
	And, thus, $A \in \Omega_N^{c_1^2 \psi}\!\!\left(\frac{R}{c_1}\right)$ in the space $X_1.$
	\end{prf}
	Let  $X = X_1 \oplus X_2,$ 	and, besides, there exists a projector $P: X \to X_1$ with the norm $\norm{P} < +\infty$ and $X_2 = \Ker P.$
	Then $X^* = X_1^\bot \oplus X_2^\bot,$ where the spaces $X_1^{\bot}, X_2^{\bot}$ are the right annihilators 
		to the spaces $X_2 \subset X$ and		$X_1 \subset X$ accordingly. 
		Continuing any functional $p\in X^*_1$ with zero on the space $X_2,$ 
		we obtain a natural isomorphism $X^*_1=X^\bot_2$. Similarly, $X^*_2=X^\bot_1$.

	\begin{lemma} \label{lemma_o_summe_space}
		Let  $X = X_1 \oplus X_2,$ 
	and, besides, there exists a projector $P: X \to X_1$ with the norm $\norm{P} < +\infty$ and 
	$X_2 = \Ker P.$
	Let the set $A_1 \subset X_1.$
	Define the set $A = \left\{x\in X | Px \in A_1 \right\}.$
	Then $N(x, A) = \left(N(Px, A_1), o\right)=\{(p,o)|\ p\in N(Px,A_1)\}$ for any $x\in A.$	
	\end{lemma}
	\begin{prf}
	Fix an arbitrary vector $x=(x_1,x_2)\in A$. \\ \noindent
	1)  Let $p = (p_1, p_2) \in N(x, A).$ Let us prove that $p_2 = o.$
	Suppose that $p_2 \neq o.$ Then there exists a non-zero vector $v=(o, v_2)$ such that
	$\bra p, v\ket = \bra p_2, v_2\ket = k\neq 0.$ According to the definition of set $A$ for any $\la \in \R$ the vector 
	$x+ \la v$ lies in $A.$ But then $\bra p, \la v\ket = \la k$ is a linear in $\la$ function and, therefore, $p \notin N(x,A).$
	Contradiction. \\ \noindent
	2) Let $p = (p_1, o) \in N(x, A).$ Let us show that $p_1 \in N(Px, A_1)= N(x_1, A_1).$
	Suppose the contrary, that $p_1 \notin N(x_1, A_1).$ Then there exists $\e > 0$ such that for any $\delta > 0$ there exists a vector $y^{\delta}_1 \in A_1$ such that $\displaystyle\norm{y_1^\delta - x_1} \men \delta$ and also $\displaystyle \bra p_1, y_1^\delta - x_1\ket \bol \e \norm{y_1^\delta - x_1}.$
	Define $y_\delta = (y_1^\delta, x_2).$ Clearly, $y_\delta \in A.$ Then from the last inequality we have that
	\begin{equation*}
	\bra p, y_\delta - x \ket = \bra (p_1, o), (y_1^\delta -x_1, o) \ket = \bra p_1, y_1^\delta - x_1\ket \bol  \e \norm{y_1^\delta - x_1} = \e \norm{y_\delta - x}.
	\end{equation*}
	Therefore, $p \notin N(x, A).$ Contradiction. \\ \noindent
	3) Let  $p_1 \in N(x_1, A_1).$ 
	Let us show that the inclusion $p=(p_1, o) \in N(x, A)$ holds.
	Fix $\e >0.$ According to the definition of the normal cone, for  $\displaystyle\frac{\e}{\norm{P}}$ there exists  $\delta > 0$ such that
for any vector $d_1 \in \B_\delta (x_1) \cap A_1$ the inequality $\bra p_1, d_1 - x_1\ket \men \frac{\e}{\norm{P}} \norm{d_1-x_1}$ holds.
	For any $y=(y_1,y_2) \in \B_{\frac{\delta}{\norm{P}}}(x) \cap A$ 
	we have $\norm{P(y-x)} = \norm{y_1 - x_1}_1 \men \norm{P} \norm{y-x} \men \delta.$ This and the previous inequality imply that
 	\begin{equation*}
	\bra p, y- x\ket = \bra p_1, y_1 - x_1\ket  \men \frac{\e}{\norm{P}} \norm{y_1 - x_1} \men \e \norm{y-x}.
	\end{equation*}
	\end{prf}
	\begin{lemma}\label{lemma_o_perenorm4}
	If, additionally, in the assertion of Lemma \ref{lemma_o_summe_space} the set $A_1$ is a complement to a convex set in  $X_1$ and 
		$A_1 \in \Omega_N^{\psi}(R)$ in $X_1,$ then $A \in \Omega_N^{\norm{P}\psi}\!\!\left(\frac{R}{\norm{P}}\right)$ in the space $X.$ 
	\end{lemma}
	\begin{prf}
	Through $\norm{\cdot}_1$ denote the norm of space $X_1^*.$ Clearly, for $p = (p_1, o) \in X^*$ the following inequality
	\begin{equation}
	\norm{p_1}_1 \men \norm{p}
	\label{eq:lemma_o_prodlenii}
	\end{equation}
	holds.
	Let $x=(x_1, x_2), y=(y_1,y_2) \in A,$ $p^x \in N(x, A) \cap \SSS,$ $p^y \in N(y, A) \cap \SSS$ and
	the following inequality holds $\norm{x_1 - y_1} \men \e,$ where $\e$ is taken from definitions of the class
	$\Omega_N^{\psi}(R).$
	Lemma \ref{lemma_o_summe_space} implies that $p^x = (p^x_1, o),$ $p^y = (p^y_1, o)$ and 
	$p^x_1 \in N(x_1, A_1),$ $p^y_1 \in N(y_1, A_1)$ in the space $X_1.$
	As $A_1$ is a complement to a convex set in the space $X_1,$ then the set $A$ is a complement to a convex set in $X.$
	Then $\bra p^x, x - y\ket \men 0.$ From this and the inequalities \reff{eq:lemma_o_prodlenii} we get that
	\begin{equation*}
	\bra p^x, x - y\ket = \bra p_1^x, x_1 - y_1 \ket \bol  \bra \frac{p_1^x}{\norm{p_1^x}_1}, x_1 - y_1 \ket.
	\end{equation*}
 	Similarly, $\bra p^y, y - x\ket  \bol  \bra \frac{p_1^y}{\norm{p_1^y}}, y_1 - x_1 \ket.$
	Summing the last two inequalities, we obtain that
	\begin{gather*}
	\bra p^x -p^y, x - y\ket \bol  \bra \frac{p_1^x}{\norm{p_1^x}_1} - \frac{p_1^y}{\norm{p_1^y}}, x_1 - y_1\ket 
	\bol - R \psif{\frac{\norm{x_1 - y_1}}{R}} \bol \\\bol  - R \psif{\frac{\norm{P}\norm{x - y}}{R}}.
	\end{gather*}
	\end{prf}

In the next lemma there are constructed examples of sets that belong to the class  $\Omega_N^{k_2\psi}(R)$ 
and do not belong to the class $\Omega_{N}(R).$
To prove it we will need some additional results from the geometry of convex sets and Banach spaces.
\begin{definition}
	For a convex set  $K \subset \R^n$ {\it the John ellipsoid} is called the maximum volume ellipsoid that is contained within $K.$
\end{definition}
It is known that for any convex  $K \subset \R^n$ such an ellipsoid exists and is unique.
Moreover, the following theorem holds (see \cite{Gruber}).
\begin{nfact}
Let $B_E$ be the John ellipsoid of the unit ball $B_n$ in the space $X_n,$ $\dim X_n = n.$
Then the following inclusions
\begin{equation}
B_E \subset B_n \subset \sqrt{n} B_E
\label{eq:ellips_John}
\end{equation} 
hold.
\end{nfact}

According to paper \cite{Figiel_1}, the convexity modulus of an arbitrary Banach space  $X$ satisfies the relation
\begin{equation}
\frac{\mcox{s}}{s^2} \men L \frac{\mcox{t}}{t^2} \quad \forall\; 0 < s \men t \men  2
\label{eq:Orlicz_Cond_mcox}
\end{equation}
for some constant $ 0< L < 4.$
\begin{lemma}\label{nadoelo}
	Let there be given a function $\psif{\cdot} \in \mathfrak{M}$ such that 
	$t^2 = \operatorname{O}(\psif{t})$ as $t \to 0.$
	Let the convexity modulus $\mcox{\cdot}$ in the Banach space $X$  satisfy the equality
	 $\mcox{t} = \omal{\psif{t}}$ as $t \to 0.$
	Then there does not exist a constant $k_2 > 0$ such that the inclusion 
	$\Omega_N^{k_2\psi}(R) \subset \Omega_{N}(R)$ holds for any $R >0$.
\end{lemma}
\begin{prf}
Suppose the contrary, that there exists a constant $k_2>0$ such that $\Omega_N^{k_2\psi}(R)\subset\Omega_N(R).$
As the condition of Lemma \ref{nadoelo} implies that $\mcox{\e} = \omal{\psif{\e}}$ as $\e \to 0,$ 
 it is sufficient to consider two cases of relations between the functions $t^2, \mcox{t} , \psif{t}$ at zero:
\begin{enumerate}
\item[1)] $t^2 = \omal{\psif{t}}$ as $t \to 0$;
\item[2)] $\mcox{t}= \omal{t^2}$ as $t \to 0.$
\end{enumerate}
Indeed, if $t^2 = \operatorname{O}(\psif{t})$ and $t^2 \neq\omal{\psif{t}}$ as $t \to 0,$
 then $\liminf\limits_{t\to 0}\frac{\psi(t)}{t^2}= c \in (0, +\infty).$
This and the equality $\mcox{t} = \omal{\psif{t}}$ as $t \to 0$ imply that $\liminf\limits_{t\to 0}\frac{\mcox{t}}{t^2}= 0.$
From this and inequality \reff{eq:Orlicz_Cond_mcox} we obtain that  $\mcox{t}= \omal{t^2}$ as $t \to 0.$

We now pass to the consideration of the two cases, mentioned above.
\\ \noindent
{\bf The case $t^2 = \omal{\psif{t}}$ as $t \to 0$.}\\ \noindent
Lemma \ref{lemma_o_perenorm} and Theorem \ref{NTh_o_mgl} imply that in a two-dimensional linear space  $Y$ there exists a norm $\norm{\cdot}_s$ such that the smoothness modulus $\mgl{X_s}{\cdot}$ of the space $X_s = (Y, \norm{\cdot}_s)$ 
satisfies the conditions $t^2 = \omal{\mgl{X_s}{t}}$ and $\mgl{X_s}{t} = \omal{\psif{t}}$ as $t \to 0.$

Therefore, there exists a constant $\gamma_0 $ such that for any $\gamma \in  (0, \gamma_0]$ the inequality
\begin{equation}
\frac{\mgl{X_s}{\frac{\gamma}{4}}}{	16\gamma} \bol 51 \gamma
\label{eq:thmain2_1743}
\end{equation}
holds.

Fix $0< \e < \min\left\{\frac{1}{2},\gamma_0\right\}.$
By the definition of the modulus of supporting convexity and \reff{suppconvprop1}, in space $X$ there exist vectors
$u, v \in\SS$  such that $v \prp u$ and
\begin{equation}
\la_X(u,v, \e) \men 2 (1 - \sqrt{1 - \e^2}) < 2\e^2.
\label{eq:thmain2_1}
\end{equation}
Let $p \in -J_1(u)$ be such that $\bra p, v\ket = 0.$
   
Consider the plane $X_2 = \lin{u,v}.$
Note that the restriction of the functional $p$ on $X_2$ does not change its norm. Denote this restriction by $p_1.$
Through $B_2$ denote the unit ball in this space with the norm induced by the space  $X.$
Consider the John ellipsoid (ellipse)  ${B_E}$ for $\frac{B_2}{\sqrt{2}}.$ 
Denote the norm, generated by the set $B_E$ as a unit ball through $\norm{\cdot}_E.$
Inclusion \reff{eq:ellips_John} implies that for any $x \in X_2$ the inequality
$2\norm{x} \bol   \norm{x}_E \bol \sqrt{2} \norm{x}$ holds.
There exists an affine transformation that transforms the John ellipsoid (ellipse) of the space $X_s$ into the ellipsoid $B_E$. 
Replacing the space $X_s$ by the isometrically isomorphic to it space, which is obtained from the affine transformation,  we consider that
the set  $B_E$ is the John ellipsoid of the unit ball  $B_s$ in the space $X_s.$
Moreover, the inequality \reff{eq:thmain2_1743} holds.
Therefore, the norms $\norm{\cdot}_E, \norm{\cdot}, \norm{\cdot}_s$ in a two-dimensional linear space are related by 
\begin{equation}
\norm{x} \men \frac{\norm{x}_E}{\sqrt{2}} \men \norm{x}_s \men  \norm{x}_E \men 2\norm{x}.
\label{eq:thmain2_2}
\end{equation}

Denote $\displaystyle{\gamma =\frac{\e}{25}}.$
In the space $X_s$, according to the definition of $\Gamma(X_s\setminus \vn B_s,\gamma, X_s)$  and the inequality \reff{eq:lemma_hypo_ball_1-1},
there exist  $d$ and $f$ from the unit ball in $X_s$ such that $\norm{d-f}_s = \gamma$ and for $q \in -J_1(f)$  the inequality
\begin{equation*}
\bra q, d -f\ket \bol \frac{1}{2} \cdot \frac{1}{2} \Gamma(X_s\setminus \vn B_s,\gamma, X_s) \bol \frac{1}{4} \mgl{X_s}{\frac{\gamma}{4}}
\end{equation*}
holds.

Executing the affine transformation $K$, which, if necessary, consists in the symmetry with respect to the main axis of the ellipse $B_E$ and rotation by some angle $\phi_0$
with respect to the ellipse $B_E$ (i.e. $K$ is such an affine transformation of the plane $X_2$ that the ellipse $B_E$ transforms into itself),
we can achieve that

\begin{enumerate}
	\item The line $l,$ supporting to $K(B_s)$ at point $Kf$, is parallel to $\Ker p = \{\la v | \;\la \in \R\}.$
	\item The projection $Kd$ on $l$ along $u$ lies in ray $\{Kf+ \la v|\; \la \bol 0\}.$ 
\end{enumerate}

 Denote $A_1 = X_s\setminus \vn K(B_s),$ $f_1= Kf,$ $d_1 = Kd.$
Using the fact that the linear functional on a two-dimensional space is defined up to a multiplicative constant by the line that is the kernel of the functional,
the functional  $q_1\in-J_1(f_1)$ (unique due to the smoothness of $B_s$)  is codirectional with the functional $p_1$  (see Fig. \ref{fig_lemma_nadoelo_1}).  
\begin{figure}[h]%
\center{
\psfrag{o}[1]{\raisebox{0.7ex}{$o$}}
\psfrag{o1}[1]{\raisebox{0.7ex}{\hspace{0.5em}$o$}}
\psfrag{f}[1]{\raisebox{2.1ex}{$f$}}
\psfrag{d}[1]{\raisebox{0.2ex}{$d$}}
\psfrag{f1}[1]{\raisebox{-1ex}{\hspace{0.1em}$f_1$}}
\psfrag{d1}[1]{\raisebox{1.4ex}{\hspace{0.75em}$d_1$}}
\psfrag{h1}[1]{\raisebox{-0.3ex}{\hspace{0.24em}$h_1$}}
\psfrag{l1}[1]{\raisebox{2.5ex}{\hspace{-1.4em}$\ell_1$}}
\psfrag{l3}[1]{\raisebox{2.1ex}{\hspace{-1.4em}$\ell_1$}}
\psfrag{l4}[1]{\raisebox{3.2ex}{$l$}}
\psfrag{l2}[1]{\hspace{0.5em}$l'$}
\psfrag{k}[1]{\raisebox{3ex}{\LARGE$K$}}
\psfrag{u}[1]{\raisebox{1.1ex}{\hspace{0.74em}$u$}}
\psfrag{u1}[1]{\raisebox{1.2ex}{\hspace{0.74em}$u$}}
\includegraphics[scale=0.7]{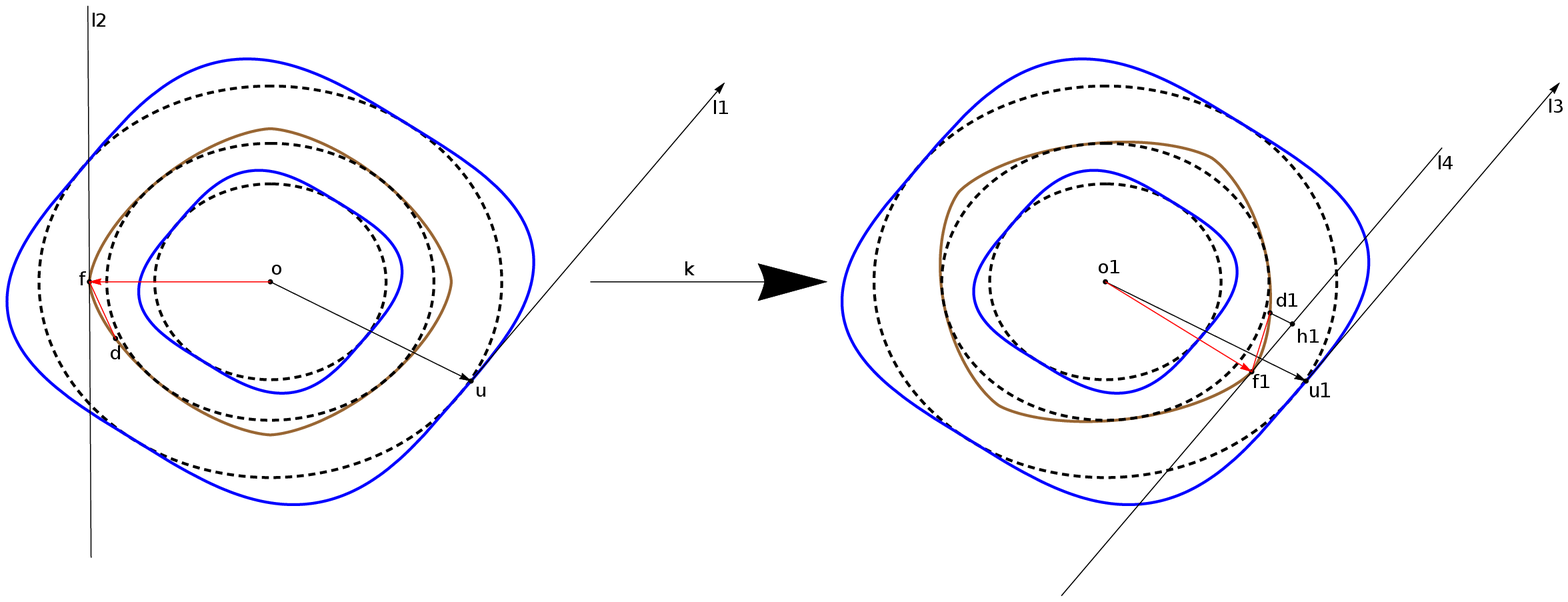}
}
\caption{ To the proof of Lemma \ref{nadoelo}.
In the Lemma notations:\\ \noindent 
-- with a black dotted line the boundaries of the ellipses = $\frac{B_E}{\sqrt{2}}, B_E, \sqrt{2}B_E$ are shown;\\ \noindent
-- with a blue line the boundaries of the balls $B_2, \frac{B_2}{2}$ are depicted; \\ \noindent
-- with brown the boundaries of balls  $B_s, K B_s$ are depicted;\\ \noindent
-- $\ell_1$ is the ray $\{u + \la v| \la \bol 0\}$ (part of the supporting line $\Ker p_1 + u$);\\ \noindent
-- $l'$ is the preimage of line $l$ (the supporting line to the ball $B_s$ at point $f$).}
\label{fig_lemma_nadoelo_1}
\end{figure}

Let $q_1 = kp_1.$ Then, as the transformation $K$ is affine, we obtain that
\begin{equation}
\bra kp_1, d_1 -f_1\ket \bol \frac{1}{4} \mgl{X_s}{\frac{\gamma}{4}}.
\label{eq:thmain2_2156}
\end{equation} 
On the other hand, it is clear that the norm $\norm{\cdot}_s^{'}$ of the space $X_s^{'}$, generated by the set  $K(B_s)$ (i.e. we consider $K(B_s)$ to be the unit ball), satisfies the following series of inequalities
\begin{equation}
\norm{x} \men \frac{\norm{x}_E}{\sqrt{2}} \men \norm{x}_s^{'} \men  \norm{x}_E \men 2\norm{x}.
\label{eq:thmain2_2157}
\end{equation}
This, particularly, implies that $1 \men k \men 2.$
As the right side of the inequality \reff{eq:thmain2_2156} is positive, we obtain that
\begin{equation}
\bra p_1, d_1 -f_1\ket \bol \frac{1}{8} \mgl{X_s^{'}}{\frac{\gamma}{4}}.
\label{eq:thmain2_2225}
\end{equation}

According to Lemma \ref{lemma_o_proektore}, in $X$ there exists a projector $P: X \to X_2$ with the norm not greater than 3,  that projects $\Ker p$ on $\Ker p_1.$
Define the set $A = \left\{x \in X| \; Px \in A_1 \right\}.$
According to  Corollary \ref{corollary_o_perenorm}, $A_1 \in \Omega_N^{\frac{1}{17}\mgl{X_s^{'}}{\cdot}}(1)$  in the space $X_s^{'}.$
Then Lemma \ref{lemma_o_perenorm_2} and the series of inequalities \reff{eq:thmain2_2157} imply that
$A_1 \in \Omega_N^{\frac{4}{17}\mgl{X_s^{'}}{\cdot}}\!\!\left(\frac{1}{2}\right)$ in the space $X_2.$
This, according to Lemma \ref{lemma_o_perenorm4}, implies that the inclusion
$A \in \Omega_N^{\frac{12}{17}\mgl{X_s^{'}}{\cdot}}\!\!\left(\frac{1}{6}\right)$ in space $X$ holds.
But by construction we have that $\mgl{X_s^{'}}{t} = \mgl{X_s}{t} = \omal{\psif{t}}$ as $t\to 0.$
Therefore, for sufficiently small $t$ the inequality $\frac{12}{17}\mgl{X_s^{'}}{6t}  \men k_2 \psif{6t}$ holds,
which implies that $A \in \Omega_N^{k_2\psif{\cdot}}\!\!\left(\frac{1}{6}\right).$
According to the assumption made, we obtain that $A \in \Omega_N\!\left(\frac{1}{6}\right).$

By construction for the point $f_1 \in X_2 \subset X$  the equality $p = (p_1, o)$ holds,
and, therefore, $p = N(f_1, A) \cap \SSS.$
Let $D = \B_{\frac{1}{6}}\!\!\left(f_1 - \frac{1}{6} u \right).$ 
According to the definition of the class $\Omega_N\!\left(\frac{1}{6}\right)$,  the following relation is satisfied
\begin{equation}
\vn D \cap A = \emptyset.
\label{eq:thmain2_2247}
\end{equation}
The inequality \reff{eq:thmain2_1} implies that $\frac{\la_X(u,v, \e)}{6} < \frac{\e^2}{3}< \frac{1}{6}.$
By the definition of $\la_X(u,v, \e)$ and the choice of vectors  $u$ и $v$  we have that
$z_1 = f_1 + \frac{\e}{6}v - \frac{\e^2}{3} u\in \vn D.$ 
In the space $X_2$ consider the projection $h_1$ of the point $d_1$ on the line $l.$
By the uniform convexity and uniform smoothness  of the space $X_2$ the segment $h_1d_1$ is parallel to $u.$
On the other hand, as the vector $d_1 \in X $  lies in $A,$ then from \reff{eq:thmain2_2247} we obtain that $d_1 \notin \vn D.$
Therefore, the segment $f_1d_1$ intersects the sphere $\partial D$ in two points.
Then, by construction, the following series of inequalities holds
$\norm{h_1 - f_1} \men\norm{f_1 - d_1} +\norm{d_1 - h_1} \men 2\norm{f_1 - d_1} \men 2\norm{f_1 - d_1}_s^{'}  = 2\gamma < \frac{\e}{6} < \frac{1}{6}.$
From this, \reff{suppconvprop2} and considering that the segment
 $h_1d_1\cap \vn D = \emptyset,$  $h_1d_1 \parallel u,$ the point $h_1$ lies on the supporting line  to the ball $D$ of radius $\frac{1}{6}$ at the point $f_1$,
we get that
\begin{equation*}
\norm{d_1-h_1} \men \frac{1}{6}\la_{X_2}\!\!\left(f_1, \frac{h_1 - f_1}{\norm{h_1 - f_1}}, 6\norm{h_1-f_1}\right) \men \norm{h_1-f_1},
\end{equation*}
and, thus, $\norm{f_1 - d_1} \men 2 \norm{h_1 - f_1}.$ 
From this,  the choice of $f_1, d_1$ and the inequality  \reff{eq:thmain2_2157}, we obtain that
\begin{equation}
0< \frac{\gamma}{2} = \frac{\norm{f_1 - d_1}_s^{'}}{2} \men \norm{f_1 - d_1} \men 2\norm{h_1 - f_1}\men 4\gamma < \frac{\e}{6}.
\label{eq:thmain2_2310}
\end{equation}
But $d_1 = f_1 + \norm{h_1-f_1}v - \bra p_1, d_1 - f_1 \ket u = f_1 + \norm{h_1-f_1}v - \bra p, d_1 - f_1 \ket u.$
The inequalities \reff{eq:thmain2_2310}, \reff{eq:thmain2_2225}, \reff{eq:thmain2_1743}  and the choice of constant $\gamma$ imply that
\begin{equation}
	\frac{\bra p, d_1 - f_1 \ket}{\norm{h_1-f_1}} \bol \frac{\mgl{X_s}{\frac{\gamma}{4}}}{16\gamma} \bol 51 \gamma > 2\e.
\label{eq:thmain2_2327}
\end{equation}
From $\norm{f_1 - h_1} < \frac{\e}{6}$ and inequality \reff{eq:thmain2_2327} we obtain that
the segments $d_1 h_1 \subset A$ and $f_1z_1 \subset D$ intersect. Denote by $y_1$ their intersection point.
Clearly, $y_1 \in \vn D \cap A$ (see Fig. \ref{fig_lemma_nadoelo_2}). This contradicts  relation  \reff{eq:thmain2_2247}.
\begin{figure}[h]%
\center{
\psfrag{o}[1]{\raisebox{1.1ex}{\hspace{0.7em}$o'$}}
\psfrag{f1}[1]{\raisebox{3.1ex}{$f_1$}}
\psfrag{d1}[1]{\raisebox{1.2ex}{\hspace{1.45em}$d_1$}}
\psfrag{h1}[1]{\raisebox{3.1ex}{\hspace{0.24em}$h_1$}}
\psfrag{y1}[1]{\raisebox{1.0ex}{\hspace{1.24em}$y_1$}}
\psfrag{v1}[1]{\raisebox{3.3ex}{\hspace{0.38em}$b$}}
\psfrag{v2}[1]{\raisebox{3.1ex}{\hspace{0.24em}$a_1$}}
\psfrag{g1}[1]{\raisebox{3.1ex}{\hspace{0.24em}$g_1$}}
\psfrag{z1}[1]{\raisebox{0.8ex}{\hspace{1.5em}$z_1$}}
\psfrag{l}[1]{\raisebox{2.1ex}{$l$}}
\includegraphics[scale=1.3]{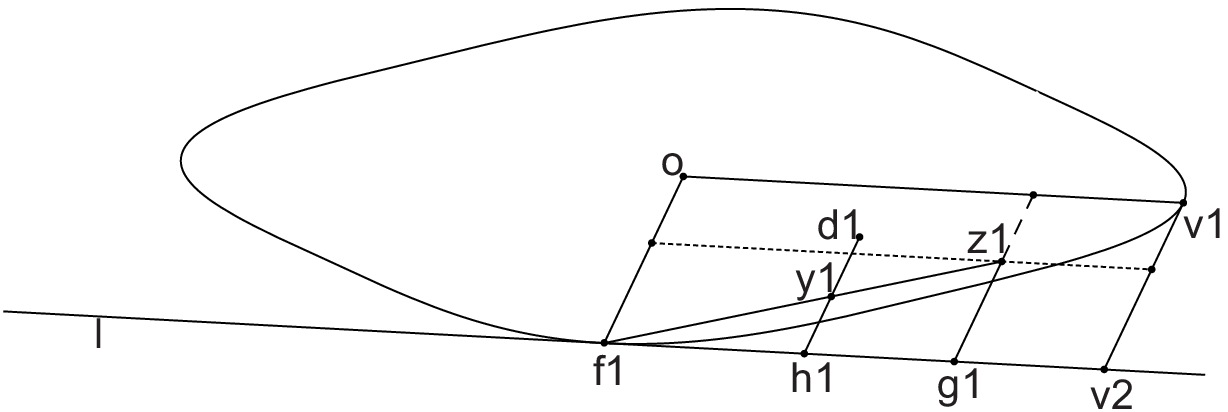}
}
\caption{ To the proof of Lemma \ref{nadoelo}.
In the notations of the lemma:\\ \noindent 
-- $o' = f_1 - \frac{1}{6}u,$ \quad $b = o' + \frac{1}{6}v,$ \quad $a_1 = f_1 + \frac{1}{6}v,$ \quad $g_1 = f_1 + \frac{\e}{6}v;$\\ \noindent
-- as $\frac{\la_X(u,v, \e)}{6} < \norm{z_1 - g_1}< \frac{1}{6},$ then $z_1$ lies inside the parallelogram $o'ba_1f_1$ and inside $D;$\\ \noindent
-- similarity implies that $\frac{\norm{z_1 - g_1}}{\norm{f_1 - g_1}} =\frac{\norm{y_1 - h_1}}{\norm{f_1 - h_1}} < \frac{\norm{d_1 - h_1}}{\norm{f_1 - h_1}}.$
}
\label{fig_lemma_nadoelo_2}
\end{figure} 

{\bf The case $\mcox{t}= \omal{t^2}$ as $t \to 0.$}\\ \noindent 
The consideration of this case does not differ a lot from the consideration of the previous one. The difference consists in the specific choice of constants. All the notations coincide with the notations in the previous case, if not mentioned otherwise.
Let $k_3 > 0$ be such that $\psif{t} \bol k_3 t^2$ in some neighborhood of zero.
Note that if $t \in [0,1]$ the inequalities
\begin{equation}
	\frac{t^2}{4} \men \mgl{H}{t} = \sqrt{1+t^2} -1 \men t^2
\label{eq:thmain2_946}
\end{equation}
hold.

Let $R_1 = \frac{17}{72} k_2 k_3.$
By assumption we have that $\mcox{t} = \omal{t^2}$ as $t \to 0.$
Fix $\gamma_0 \in (0, 1)$ such that for any $t \in (0, \gamma_0)$ the inequality
$\mcox{t} <  \frac{1}{3}\cdot\frac{R_1}{2^{12}}t^2$ holds.

By the definition of the modulus of supporting convexity and \reff{suppconvprop1} for
$$\displaystyle{\e \in \left(0, \min\!\left\{\frac{\gamma_0}{2}, \frac{1}{\sqrt{R_1}}\right\}\right)}$$ 
there exist
$u, v \in\SS$ in the space $X$ such that $v \prp u$ and
\begin{equation}\label{eq:thmain2_111}
\la_X(u,v, \e) \men 2 \mcox{2\e} < \frac{1}{3}\cdot\frac{R_1}{512}\e^2 < 1.
\end{equation}
The space $X_2$ and the functionals  $p, p_1$ are defined similarly to the previous case.
We consider the unit ball in $X_s^{'}$ to be the set $B_E.$
Fix $\gamma = \frac{R_1\e}{3}.$
Similarly to the previous case, we take the points $f_1, d_1, h_1,$ and construct the set $A.$
As in the previous case, $A \in \Omega_N^{\frac{12}{17}\mgl{X_s^{'}}{\cdot}}\!\!\left(\frac{1}{6}\right).$
As the ball in $X_s^{'}$ is Euclidean and considering the left inequality in  series  \reff{eq:thmain2_946} we have that
\begin{equation}\label{eq:thmain2_1003}
\bra p_1, d_1 -f_1\ket \bol \frac{1}{4} \mgl{X_s^{'}}{\frac{\gamma}{4}} \bol \frac{1}{256}\gamma^2.
\end{equation}

From the fact that the norm in  $X_s^{'}$ is Euclidean, the right hand side of inequality  \reff{eq:thmain2_946} and the choice of constants $k_3$ и $R_1$ 
we obtain that
\begin{equation*}
R_1k_2 \psif{\frac{t}{R_1}} \bol \frac{k_2 k_3 t^2}{R_1} = \frac{72 t^2}{17} \bol \frac{1}{6} \cdot\frac{12}{17} \mgl{X_{s'}}{6t}.
\end{equation*}
Thus, $A \in \Omega_N^{k_2\psif{\cdot}}\!\!\left(R_1\right).$
The assumption  implies that $A \in \Omega_N(R_1).$

Let $D = \B_{R_1}(f_1 -R_1 u).$ 
Inequality \reff{eq:thmain2_111} and the choice of constant $\e$ imply that $R_1\la_X(u,v, \e) <  \frac{1}{3}\cdot\frac{R_1^2}{512}\e^2 < R_1.$
By the definition of $\la_X(u,v, \e)$ and the choice of vectors $u$ and $v$  we get that
$z_1 = f_1 + R_1\e v - \frac{1}{3}\cdot\frac{R_1^2}{512}\e^2 u \in \vn D.$ 

Similarly to the previous case, the estimations of the norm  $h_1 - f_1:$ are made
\begin{equation}\label{eq:thmain2_1006}
\frac{\gamma}{4}  \men \norm{h_1 - f_1}\men 2\gamma < R_1\e.
\end{equation}
Furthermore, $d_1 =  f_1 + \norm{h_1-f_1}v - \bra p, d_1 - f_1 \ket u.$
Inequalities \reff{eq:thmain2_1006}, \reff{eq:thmain2_1003} and the choice of constant $\gamma$ 
imply that
\begin{equation}
	\frac{\bra p, d_1 - f_1 \ket}{\norm{h_1-f_1}} \bol   \frac{1}{512} \gamma =  \frac{1}{3}\cdot\frac{1}{512} R_1\e.
\label{eq:thmain2_232711}
\end{equation}
By the inequality  $\norm{f_1 - h_1} < R_1 \e$ and  inequality \reff{eq:thmain2_232711} we get that
the segments $d_1 h_1 \subset A$ and $f_1z_1 \subset D$ intersect. Denote their intersection point by $y_1.$
Clearly, $y_1 \in \vn D \cap A.$ This contradicts the relation \reff{eq:thmain2_2247}.
\end{prf}

\noindent {\bf Proof of the main theorem.}\\ \noindent
	Let us show that the proof 1 $\Rightarrow$ 2 in Theorem \ref{Th_npsi_to_prox} follows from Lemma \ref{nadoelo}.

Inequalities \reff{eq:Orlicz_Cond_mcox}, \reff{eq:Orlicz_Cond} imply that if for the function 
$\psi \in \mathfrak{M}_2$ the relation
$\psi(\e) \neq\operatorname{O}(\mcox{\e})$ holds as $\e \to 0,$ then 
$\mcox{\e} = \omal{\psif{\e}}$ as $\e \to 0.$ 
\bbox

\subsection{Proof of Theorem \ref{mainTh_3}}
\begin{lemma} \label{lemma_iskl_sluch_mainTh}%
Let $X$ be a uniformly smooth and uniformly convex Banach space.
Let there be given a function $\psif{\cdot} \in \mathfrak{M}$ %
such that $t^2 = \operatorname{O}(\psif{t})$ as $t \to 0.$%
If there exists a constant $k_2 > 0 $ such that the inclusion $\Omega_N^{k_2\psi}(R) \subset \Omega_{N}(R)$ holds,
then $\mcox{\e} \asymp \e^2$ at zero and  $\liminf\limits_{t\to 0}\frac{\psi(t)}{t^2}= c \in (0, +\infty).$
\end{lemma}%
\begin{prf}
Applying Lemma \ref{nadoelo}, we get that  $\mcox{\e} \neq \omal{\psif{\e}}$ as $\e \to 0.$
But for the function  $\psif{\cdot} \in \mathfrak{M}$  such that  $t^2 = \operatorname{O}(\psif{t})$ as $t \to 0,$
the relation $\mcox{\e} \neq \omal{\psif{\e}}$ as $\e \to 0$ may hold if and only if  $\mcox{\e} \asymp \e^2$ as $\e \to 0$ and 
$\liminf\limits_{t \to 0}\frac{\psif{t}}{t^2} = c \in (0, +\infty).$
\end{prf}
\noindent 
{\bf Proof of Theorem \ref{mainTh_3}.}\\ \noindent 
Applying Theorem \ref{Th_prox_to_npsi}, we obtain that in some neighborhood of zero for some constant  $k > 0$ 
the inequality
\begin{equation}
\label{thelastbutnotleast}
\psif{\e} \bol k \mglx{\e}
\end{equation}
holds.
This and the Day-Nordlander theorem (see \cite{DiestelEng})  imply that $\e^2 = \operatorname{O}(\psif{\e})$ as $\e \to 0.$ 
Now, applying Lemma \ref{lemma_iskl_sluch_mainTh}, we obtain that $\mcox{\e} \asymp \e^2$ as $\e \to 0$ and
$\liminf\limits_{t \to 0}\frac{\psif{t}}{t^2} = c \in (0, +\infty).$
This, inequality \reff{thelastbutnotleast} and the fact that the smoothness modulus satisfies the Figiel condition 
(see Remark \ref{rem_figiel}) imply that $\mglx{\e} \asymp \e^2$ as $\e \to 0.$

Since $\mcox{\e} \asymp \mglx{\e} \asymp \e^2$ as $\e \to 0,$  $X$ is isomorphic to a Hilbert space  (see \cite{Borwein_book1},\,\S 5.5.6).
\bbox

\end{document}